 \documentclass[final,3p,times]{elsarticle}

\usepackage{graphicx}
\usepackage{amsmath,amssymb}
\usepackage{subeqn}
\usepackage{booktabs}
\usepackage{subfigure}

\def\d{\mathrm{d}} \def\dd{\partial}
\newcommand{\spec}[1]{\mathrm{#1}}

\newtheorem{dfn}{Definition}

\biboptions{comma,square,sort&compress}

\journal{Journal of Computational Physics}

\begin{document}

\begin{frontmatter}



\title{Minimal Curvature Trajectories: Riemannian Geometry Concepts for Model
  Reduction in Chemical Kinetics}

\author[label1,label2]{Dirk Lebiedz\corref{cor1}}
\ead{dirk.lebiedz@biologie.uni-freiburg.de}
\author[label2,label3]{Volkmar Reinhardt}
\author[label2]{Jochen Siehr}

\cortext[cor1]{Corresponding author}

\address[label1]{Center for Analysis of Biological Systems (ZBSA), University of Freiburg, Habsburgerstra\ss{}e 49, 79104 Freiburg, Germany}
\address[label2]{Interdisciplinary Center for Scientific Computing (IWR), University of Heidelberg, Im Neuenheimer Feld 368, 69120 Heidelberg, Germany}
\address[label3]{SEW-EURODRIVE GmbH \& Co KG, Ernst-Blickle-Str.~42, 76646 Bruchsal, Germany}

\begin{abstract}
  In dissipative ordinary differential equation systems different time scales
  cause anisotropic phase volume contraction along solution
  trajectories. Model reduction methods exploit this for simplifying chemical
  kinetics via a time scale separation into fast and slow modes.  The aim is
  to approximate the system dynamics with a dimension-reduced model after
  eliminating the fast modes by enslaving them to the slow ones via
  computation of a slow attracting manifold. We present a novel method for
  computing approximations of such manifolds using trajectory-based
  optimization. We discuss Riemannian geometry concepts as a basis for
  suitable optimization criteria characterizing trajectories near slow
  attracting manifolds and thus provide insight into fundamental geometric
  properties of multiple time scale chemical kinetics. The optimization
  criteria correspond to a suitable mathematical formulation of ``minimal
  relaxation'' of chemical forces along reaction trajectories under given
  constraints. We present various geometrically motivated criteria and the
  results of their application to three test case reaction mechanisms serving
  as examples. We demonstrate that accurate numerical approximations of slow
  invariant manifolds can be obtained.
\end{abstract}


\end{frontmatter}


\section{Introduction}
The need for model reduction in chemical kinetics is mainly motivated by the
fact that the computational effort for a full simulation of reactive flows,
e.g.\ of fluid transport involving multiple time scale chemical reaction
processes, is extremely high. For detailed chemical reaction mechanisms
involving a large number of chemical species and reactions, a simulation in
reasonable computing time requires reduced models of chemical kinetics
\cite{Warnatz2006}.

However, model reduction is often also of general interest for theoretical
purposes in mathematical modeling. Reduced models are intended to describe
some essential characteristics of dynamical system behavior while fading out
other issues. Therefore, they often allow a better insight into complicated
reaction pathways, e.g.\ in biochemical systems \cite{Lebiedz2005c}, and their
nonlinear dynamics.

\begin{figure}
  \begin{center}
    \includegraphics[width=8.cm]{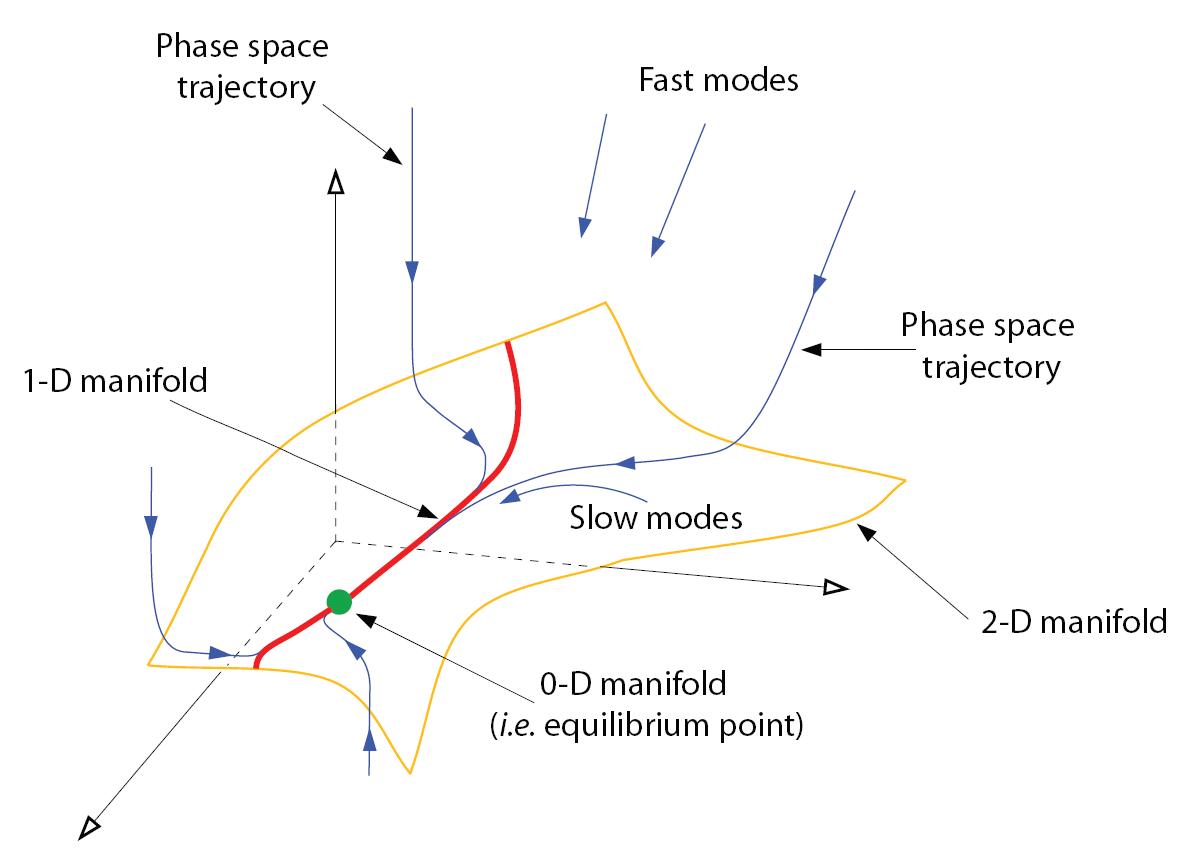}
  \end{center}
  \caption{\label{f:powers_manifold}Illustration of trajectories relaxing
    successively onto a 2-D manifold and a 1-D manifold before converging to
    equilibrium. Figure courtesy of A.N.\ Al-Khateeb, J.M.\ Powers, S.\
    Paolucci (private communication).}
\end{figure}

In dissipative ordinary differential equation systems modeling chemical
reaction kinetics different time scales cause anisotropic phase volume
contraction along solution trajectories. This leads to a bundling of
trajectories near ``manifolds of slow motion'' of successively lower dimension
as time progresses, illustrated in Figure \ref{f:powers_manifold}.
Many model reduction methods exploit this for simplifying chemical kinetics
via a time scale separation into fast and slow modes.  The aim is to
approximate the system dynamics with a dimension-reduced model after
eliminating the fast modes by enslaving them to the slow ones via computation
of slow attracting manifolds.

Very early model reduction approaches like the quasi steady-state (QSSA) and
partial equilibrium assumption (PEA) \cite{Warnatz2006} performed ``by hand'',
have set the course for modern numerical model reduction methods that
automatically compute a reduced model without need for detailed expert
knowledge of chemical kinetics by the user. Many of these modern techniques
are explicitly or implicitly based on a time-scale analysis of the underlying
ordinary differential equation (ODE) system with the purpose to identify a
slow attracting manifold in phase space where -- after a short initial
relaxation period -- the system dynamics evolve. For a comprehensive overview
see e.g.\ \cite{Gorban2005} and references therein.

In 1992 Maas and Pope introduced the ILDM-method \cite{Maas1992} which became
very popular and widely used in the reactive flow community, in particular in
combustion applications. Based on a singular perturbation approach, a local
time-scale analysis is performed on the Jacobian of the system of ordinary
differential equations modeling chemical kinetics. Fast time scales are
assumed to be fully relaxed and after a suitable coordinate transformation
fast variables are computed as a function of the slow ones by solving a
nonlinear algebraic equation system. For recent developments and extensions of
the ILDM method see e.g.\ \cite{Bykov2008} and references therein.

Another popular technique is computational singular perturbation (CSP) method
proposed by Lam in 1985 \cite{Lam1985, Lam1994}. The basic concept of this
method is a representation of the dynamical system in a set of ``ideal'' basis
vectors such that fast and slow modes are decoupled.  In contrast to ILDM and
CSP some iterative methods came into application that are not based on
explicit computation of a time scale separation, but rather an evaluation of
functional equations suitably describing the central characteristics of a slow
attracting manifold, for example invariance and stability. Examples are
Fraser's algorithm \cite{Davis1999, Fraser1988,Nguyen1989} and the method of
invariant grids \cite{Chiavazzo2007, Gorban2005a, Gorban2005}. Other widely
known and successful methods are the constrained runs algorithm
\cite{Gear2005,Zagaris2009}, the rate-controlled constrained equilibrium
(RCCE) method \cite{Keck1971}, the invariant constrained equilibrium edge
preimage curve (ICE-PIC) method \cite{Ren2005,Ren2006a}, and
flamelet-generated manifolds \cite{Delhaye2007, Oijen2000}. In
\cite{Mease2003} finite time Lyapunov exponents and vectors are analyzed for
evaluation of timescale information. Mitsos et al.\ formulated an integer
linear programming problem explicitely minimizing the number of species in the
reduced model subject to a given error constraint \cite{Mitsos2008}.

It is obvious that non-local information on phase space dynamics has to be
taken into account to get accurate approximations of slow attracting manifolds
in the general case. Reaction trajectories in phase space that are solutions
of the ODE system describing chemical kinetics and uniquely determined by
their initial values bear such information. Based on Lebiedz' idea to search
for an extremum principle that distinguishes trajectories on or near slow
attracting manifolds, we apply an optimization approach for computing such
trajectories \cite{Lebiedz2004c, Lebiedz2006b}. Various optimization criteria
have been suggested \cite{Reinhardt2008}, systematically investigated and the
trajectory-based approach has been extended to the computation of manifolds of
arbitrary dimension via parameterized families of trajectories.

This paper derives and comprehensively discusses various geometrically
motivated objective criteria for computing trajectories approximating slow
attracting manifolds in chemical kinetics as a solution of an optimization
problem. The corresponding objective functionals are supposed to implicitly
incorporate essential characteristics of slow attracting manifolds related to
a minimal remaining relaxation of chemical forces along trajectories on these
manifolds. We consider the picture of abstract ``dissipative chemical forces''
imagined to drive the single elementary reaction steps
\cite{Kondepudi1998}. Due to energy dissipation these forces successively
relax while the chemical system is approaching equilibrium.  The successive
relaxation of these forces causes curvature in the reaction trajectories (in
the sense of velocity change along the trajectory). A slow 1-D manifold in
this picture would correspond to a minimally curved reaction trajectory along
which the remaining relaxation of chemical forces is minimal while approaching
chemical equilibrium.

In particular, we propose and motivate an optimization criterion suitably
measuring curvature which is rooted in a thermodynamically motivated
Riemannian geometry specifically defined for chemical reaction kinetics and
based on the Second Law of Thermodynamics. This metric provides the phase
space of chemical reaction kinetics with a geometry specifically capturing the
structure of chemical kinetic systems \cite{Gorban2004}.

The proposed model reduction method is automatic. The user has to provide only
the desired dimension of the reduced model and the range of concentrations of
the reaction progress variables supposed to parameterize the reduced
model. For the application examples presented in this work, the numerical
optimization algorithm shows fast convergence independent of the initial
values chosen for numerical initialization. The optimization problems seem to
be convex for the example systems presented in this paper (see ``optimization
landscapes'' in Section \ref{sss:ds_ol} and \ref{sss:o3_ol}) and would then
have a unique solution corresponding to a global minimum of the objective
functional.  An advantage of the presented trajectory optimization approach
over local time scale separation methods like QSSA and ILDM is the fact that
it produces smooth manifolds and whole 1-D manifolds (trajectories) as a
solution of a single run of the optimization algorithm. Methods based on
explicit local time scale separation might yield non-smooth manifolds when the
fast-slow spectral decomposition changes its structure.  Furthermore, the
formulation as an optimization problem assures results even under conditions
where the time scale separation is small and many common model reduction
methods fail or numerical solutions are difficult to obtain.

\section{Trajectory-based optimization approach}

As described in the introduction, the key of the method presented here is the
exploitation of global phase space information contained in the behavior of
trajectories on their way towards chemical equilibrium. This information can
be used within an optimization framework for identifying suitable reaction
trajectories approximating slow invariant and attracting manifolds (SIM).  A
suitable formulation as the numerical solution of an optimization problem
assures the existence of a reduced model irrespective of assumptions on the
time scale spectrum, its structure and the dimension of the reduced model and
sophisticated optimization software can be used for the numerical solution of
the problem. The central idea behind our approach is that the optimization
criterion for the identification of suitable trajectories should represent the
assumption that chemical forces are -- under the given constraints -- already
maximally relaxed along trajectories on the slow attracting manifold. From the
opposite point of view this means that the remaining relaxation of chemical
forces along the trajectories is minimal while approaching chemical
equilibrium. This means that the velocity change caused by chemical force
relaxation is minimal which is intuitively close to the notion of a slow
manifold. Various ideas for the formulation of suitable optimization criteria
are conceivable.

Mathematically the basic problem can be formulated as
\begin{subequations} \label{eq:op}
  \begin{equation} \label{eq:op:of} \min_{c(t)} \int_0^{t_{\textrm f}}
    \Phi\left(c(t)\right) \; \d t
  \end{equation}
  \textrm{subject to}
  \begin{align}
    \frac{\d c(t)}{\d t} &= f\left(c(t)\right)  \label{eq:op:dyn}\\
    0 &= g\left(c(0)\right)  \label{eq:op:cr}\\
    c_k(0) &= c_k^0,\qquad k \in I_{\text{fixed}}. \label{eq:op:pv}
  \end{align}
\end{subequations}
The variables $c_k$ denote the concentrations of chemical species, and
$I_{\mathrm{fixed}}$ is an index set that contains the indices of variables
with fixed initial values (the so-called reaction progress variables) chosen
to parameterize the reduced model, i.e.\ the slow attracting manifold to be
computed. Thus, those variables representing the degrees of freedom within the
optimization problem are the initial value concentrations of the chemical
species $c_k(0), k\notin I_{\rm fixed}$. The process of determining $c^0_k,k
\notin I_{\rm fixed}$ from $c^0_k,k \in I_{\rm fixed}$ is known as ``species
reconstruction'' and represents a function mapping the reaction progress
variables to the full species composition by determining a point on the slow
attracting manifold. In our approach, species reconstruction is possible
locally, i.e.\ without being forced to compute the slow attracting manifold as
a whole. The system dynamics (chemical kinetics determined by the reaction
mechanism) are described by (\ref{eq:op:dyn}) and enter the optimization
problem as equality constraints. Hence an optimal solution of (\ref{eq:op})
always satisfies the system dynamics of the full ODE system. Chemical element
mass conservation relations that have to be obeyed due to the law of mass
conservation are collected in the linear function $g$ in (\ref{eq:op:cr}). The
initial concentrations of the reaction progress variables are fixed via the
equality constraint (\ref{eq:op:pv}).

When asymptotically approaching the equilibrium point $c^{\mathrm{eq}}$, which
is a stable fixed point attractor in a closed chemical system, the system
dynamics become infinitely slow and equilibrium will never be reached
exactly. By approximating the equilibrium point within a surrounding of small
radius $\varepsilon > 0 $ for the concentration of chemical species (e.g.\ the
reaction progress variables) by an additional constraint $\vert c_k(t_{\textrm
  f}) - c_k^{\mathrm{eq}} \vert \leqslant \varepsilon$ (equilibrium
composition: $c^{\mathrm{eq}}$), the free final time $t_{\textrm f}$ can be
automatically determined within the optimization problem assuring that this
additional inequality constraint is fulfilled. However, in practical
applications it is usually sufficient to choose $t_{\textrm f}$ large enough
for the final point of the integration to be close to equilibrium. The
objective functional $\Phi(c(t))$ in (\ref{eq:op:of}) characterizes the
optimization criterion which will be discussed later in detail.

The key idea of our approach to model reduction in chemical kinetics is found
in the fact that suitable trajectories can be used to span slow attracting
invariant manifolds. The approximated {SIM} can then be used as a reduced
model of the underlying {ODE} model, for example via a look-up table for
points on the slow manifold. This reduced model is parametrized by the
reaction progress variables (coordinate axes) which find a fully natural
realization in our formulation as trajectory initial concentrations
(\ref{eq:op:pv}).

\subsection{Numerical Methods: Multiple Shooting in a Parametric Optimization
  Setting} \label{ss:numerics} The optimization problem (\ref{eq:op}) can be
solved as a standard nonlinear optimization problems (NLP), for example via
the sequential quadratic programming (SQP) method \cite{Powell1978}. However,
one has to decide how to treat the differential equation constraint and the
objective functional. The easiest way is a decoupled iterative approach, a
full numerical integration of the ODE model with the current values of the
variables subject to optimization. This is called the sequential (or single
shooting) approach since it fully decouples simulation of the model and
optimization. However, it is often beneficial to have an ``all at once''
approach that couples simulation and optimization via discretization of the
ODE constraint. This simultaneous approach has the advantage of introducing
more freedom into the optimization problem since the differential equation
model does not have to be solved exactly in each iteration of the
optimization. A beneficial approach to couple the ODE constraint to the
optimization is the multiple shooting method.  Here, the time interval is
subdivided into several multiple shooting intervals and additional degrees of
freedom are introduced at the initial points of each interval. On each
multiple shooting interval an independent initial value problem is solved via
numerical integration. Additional ``matching condition''-equality constraints
at the level of the optimization problem assure continuity of the final
solution trajectory between the multiple shooting subintervals. For all
results in this paper the multiple shooting approach introduced by Bock and
Plitt \cite{Bock1987,Bock1984} is used.

The SQP method basically can be interpreted as Newton's method applied to the
Karush--Kuhn--Tucker (KKT) necessary optimality conditions of the NLP (see
e.g.\ \cite{Nocedal2006}) and requires the computation of derivatives. For the
numerical approximation of these derivatives by finite difference methods,
along with the nominal ODE solution trajectory $n$ perturbed trajectories have
to be computed, where $n$ is the dimension of the ODE system. To avoid the
dependence of the resulting derivative on the adaptive discretization schemes
of these trajectories provided by an automatic step size control in the
numerical integration routine, the perturbed trajectories are evaluated on the
same grid as the nominal trajectory. This approach is called internal
numerical differentiation (IND) \cite{Bock1981}. As the systems considered
here are chemical reaction systems which are usually stiff systems, the
integration itself is performed by DAESOL \cite{Albersmeyer2008, Bauer1997}, a
multistep backward differentiation formula (BDF) differential algebraic
equation (DAE) solver.  For all computations presented in the results section
of this paper, the software package {MUSCOD-II}
\cite{Bock1984,Leineweber2003,Leineweber2003a} has been used.

For the computation of slow attracting manifolds of dimension larger than one,
a sequence of problems of type (\ref{eq:op}) has to be solved for different
initial values of the reaction progress variables in (\ref{eq:op:pv}). For
this purpose, we use a parametric optimization framework, where neighboring
problems are efficiently initialized with the previous optimal solution.
Through this continuation method embedding the problem into a parametric
family of optimization problems, the computation of a family of optimal
trajectories spanning a higher-dimensional manifold can be significantly
accelerated. Such embedding strategy was originally developed and implemented into the
package {MUSCOD-II} by Diehl et al.\ in \cite{Diehl2005, Diehl2002b} for fast
online optimization, especially real-time optimal control. A variant of this
implementation has been used for the results presented in this paper.

\subsection{Optimization criteria}\label{ss:OptCrit}
Naturally, the choice of the criterion $\Phi(c(t))$ crucially affects both
success and degree of accuracy of the computed approximations of the slow
attracting manifold. A useful criterion $\Phi(c(t))$ should at least fulfill
the following three requirements:
\begin{enumerate}
\item $\Phi$ should be physically motivated and describe in a suitable sense
  the extent of relaxation of ``chemical forces'' or ``dynamical modes'' in
  the evolution of reaction trajectories towards equilibrium.
\item $\Phi$ should be computable from easily accessible data contained in
  standard models of chemical reaction mechanisms (e.g.\ reaction rates,
  chemical source terms and their derivatives, thermodynamic data).
\item $\Phi$ should be twice continuously differentiable along reaction
  trajectories.
\end{enumerate}
Another desirable but not necessary property, which is related to the
invariance of a manifold is the following consistency property. If a criterion
is consistent in the sense of the following definition, the manifold
computed as a solution of the optimization problem is positively invariant
which means that trajectories starting on the manifold at time $t_0$ will stay
on the manifold for all $t \geqslant t_0$.

\begin{dfn}[Consistency property]\label{def:consis}
  Suppose an optimal trajectory $\tilde{c}(t)$ has been computed as a solution
  of (\ref{eq:op}). Take the concentrations of the progress variables
  $\tilde{c}(t_1)$ at some time $t_1 > 0$ as new initial concentrations for
  (\ref{eq:op:pv}) and solve (\ref{eq:op}) again. If $\hat{c}(t) = \tilde{c}(t
  + t_1)$ holds for the resulting optimal trajectory $\hat{c}(t)$, we call the
  optimization criterion $\Phi$ consistent.
\end{dfn}

\begin{figure}
  \begin{center}
    \includegraphics[width=8.cm]{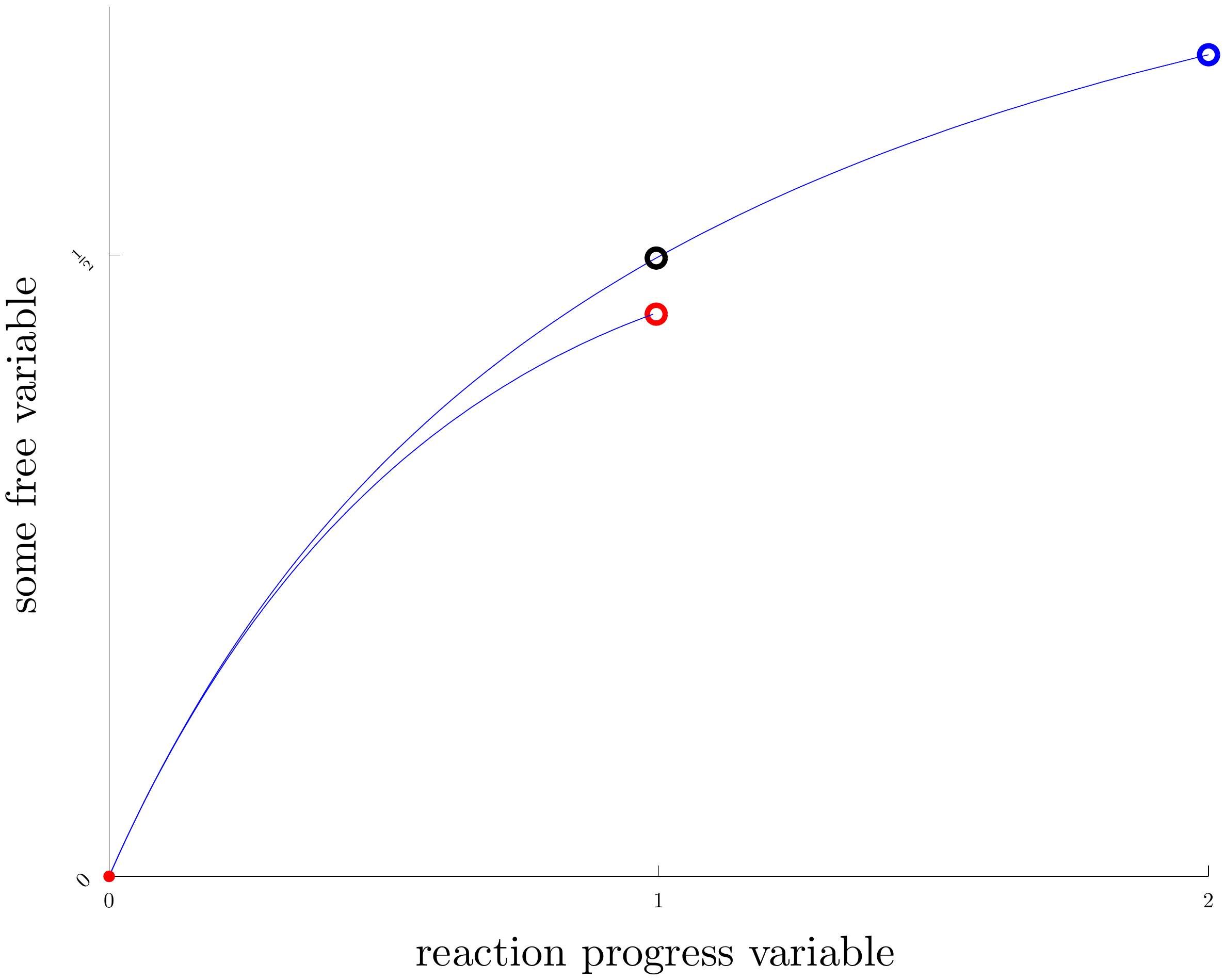}
  \end{center}
  \caption{\label{f:consis} Illustration of the consistency property: An
    optimization problem has been solved for a fixed value (=$2.0$ here) for
    the progress variable. Its solution is the trajectory $\tilde{c}(t)$
    starting from the blue circle and converging towards equilibrium (here
    the coordinate origin, red dot). At a later point in time $t_1 > 0$, the progress
    variable is fixed to the value on the trajectory $\tilde{c}_k(t)$ (=$1.0$
    here), $k \in I_{\rm fixed}$, and the optimization problem is solved
    again. If the new solution $\hat{c}(t)$ coincides with the remaining part
    of the previous one such as the trajectory starting at the black circle, we
    call the criterion (\ref{eq:op:of}) consistent, otherwise (e.g.\ the
    trajectory starting at the red circle) it is called inconsistent.}
\end{figure}

The consistency property, illustrated in Fig.\ \ref{f:consis}, can be used as
an accuracy test for the computed manifold, because the correct attracting
{SIM} should be invariant by definition.  However, it poses a strong demand
that is not {\em a priori} incorporated into the problem formulation
(\ref{eq:op}) and will not be fulfilled in general for solutions of the
optimization problem. Nevertheless, an invariant manifold can in principle be
constructed in our approach even without a consistent criterion by solving
(\ref{eq:op}) for initial values $c_k^0,\; k \in I_{\mathrm{fixed}}$ on the
boundary of a desired domain and spanning the low-dimensional manifold by the
resulting trajectories.

\subsubsection{Curvature-based Relaxation Criteria}
\label{ss:optimcrit}
As pointed out before, a suitable optimization criterion $\Phi(c(t))$ should
characterize the extent of relaxation of ``chemical forces''.  Fundamentally
rooted criteria of this type can be derived on the basis of the concept of
curvature of trajectories in phase space measured in a suitable metric. From a
physical point of view curvature is closely related to the geometric
interpretation of a force and its action on the system dynamics. This picture
has a long historical tradition.

One of the most popular examples is Einstein's general theory of relativity
\cite{Einstein1916} which proposes the idea that gravitational force is
replaced by a ``geometric picture''.  Einstein's general theory of relativity
relates the special theory of relativity and Newton's law of universal
gravitation with the insight that gravitation can be described by curvature of
space-time. Space-time is treated as a four-dimensional manifold whose
curvature is due to the presence of mass, energy, and momentum.

But even long before Einstein, the concept of curvature has already been
related to the concept of force in physics.  In 1687 Sir Isaac Newton
published the laws of motion in his work ``Philosophiae Naturalis Principia
Mathematica''. In a differential formulation Newton's second law can be stated
as $F = m \cdot a,$ where $m$ is mass, $a$ is acceleration, and $F$ is
force. Since the acceleration $a$ is the second derivative of the state
variable $x(t)$ with respect to time, $a=\ddot{x}$, and thus contains
information about the curvature of $x$; Newton's law is the first one to
directly relate force to curvature.

In this context it is important to remark that equations of motion in
classical mechanics can also be described by a variational principle,
Hamilton's principle of least action. In Lagrange-Hamilton mechanics
\cite{Goldstein1980}, the trajectory of an object is determined in such a way
that the action (which is defined as the integral of the Lagrangian over time,
where the Lagrangian is the difference of kinetic energy and potential energy)
is minimal.

A central issue in this paper is to transfer the principle of ``force =
curvature'' to the field of chemical systems and look for a corresponding
variational principle characterizing the kinetics along a slow attracting
manifold. In chemical systems dissipative forces are active. The different
time scales of dynamic modes result in an anisotropic force relaxation in
phase space. This force relaxation changes the reaction velocity. Inspired
by an analogy observation with Newton's geometric interpretation of a
force as a second derivative of a trajectory with respect to time, we regard
the second time-derivative of the chemical composition $c(t)$ characterizing
the rate of change of reaction velocity through relaxation (dissipation) of
chemical forces:

\begin{equation}
  \dot{c} = f (c), \quad \ddot{c} = \frac{\d \dot{c}}{\d t} = \frac{\d
    \dot{c}}{\d c} \cdot \frac{\d c}{\d t} = J_f(c) \cdot f.
  \label{eq:forcecurv}
\end{equation}

We consider the tangent (reaction velocity) vectors $\dot{c}(t) = f(c(t))$ of
reaction trajectories. The relaxation of chemical forces results in a change
of $\dot{c}(t)$ along a trajectory on its way towards chemical
equilibrium. This change along the trajectory may be characterized by taking
the directional derivative of the tangent vector of the curve $c(t)$ with
respect to its own direction $v\coloneqq\frac{\dot{c}}{\|\dot{c}\|_2}=\frac{f}{\|f
  \|_2}$. Mathematically this can be formulated as
\begin{equation*}
  D_v \dot{c}(t) \coloneqq \frac{\d}{\d \alpha} f(c(t) + \alpha v)
  \Big|_{\alpha=0}  = J_{f}(c) \cdot \frac{f}{\|f\|_2},
\end{equation*}
with $J_f(c)$ being the Jacobian of the right hand side $f$ evaluated at
$c(t)$ and $\|\cdot\|_2$ denoting the Euclidian norm. Hence, we may choose the
optimization criterion
\begin{equation} \label{eq:crit_c2} \Phi_{\textrm{A}}(c) = \frac{\|J_f(c) \;
    f\|_2}{\|f\|_2}
\end{equation}
in the formulation (\ref{eq:op:of}). This criterion bears some resemblance to
the recently published method of stretching-based diagnostics
\cite{Adrover2007} and its application for model reduction (SBR-method). The
authors use an expression closely related to criterion (\ref{eq:crit_c2})
which measures the stretching of vector fields in the tangent bundle of
manifolds.

The natural way for the evaluation of criterion (\ref{eq:crit_c2}) in the
formulation of the objective functional (\ref{eq:op:of}) would be a path
integral along the trajectory towards equilibrium
\begin{equation*}
  \int_{l(0)}^{l(t_{\textrm f})} \Phi(c(l(t))) \;\d l(t),
\end{equation*}
where $l(t)$ is the Euclidian length of the curve $c(t)$ at time $t$ given by
\begin{equation*}
  l(t) = \int_0^t \|\dot{c}(\tau)\|_2 \;\d \tau.
\end{equation*}
This results in the reparametrization
\begin{equation} \label{eq:repara} \d l(t) = \| \dot{c}(t) \|_2 \;\d t.
\end{equation}
The objective used in (\ref{eq:op:of}) would be (using (\ref{eq:forcecurv})):
\begin{equation} \label{eq:c2of} \int_0^{t_{\textrm f}} \|J_f(c) \; f \|_2 \;
  \d t=\int_0^{t_{\textrm f}} \|\ddot{c} \|_2 \; \d t.
\end{equation}

However, an alternative norm for the evaluation of $\|J_f(c)\;f\|$ might be
taken into account, which has already been used by Weinhold in
\cite{Weinhold1975} and is motivated from thermodynamics. This norm is also
known as Shashahani norm \cite{Shahshahani1979} and is employed for model
reduction purposes in \cite{Gorban2004}. In this norm the criterion adapted
from (\ref{eq:crit_c2}) can be written as
\begin{equation} \label{eq:crit_cw} \Phi_{\textrm{B}}(c) = \frac{\|J_f(c)\;
    f\|_W}{\|f\|_W} = \frac{(f^{\textrm{T}} J_f^{\textrm{T}}(c)\cdot
    \mathrm{diag}(1/c_i)\cdot J_f(c) f)^{1/2}}{(f^{\textrm{T}} \cdot
    \mathrm{diag}(1/c_i)\cdot f)^{1/2}}
\end{equation}
with $W=\mathrm{diag}(1/c_i)$ being the diagonal matrix with diagonal elements
$1/c_i$. This criterion brings thermodynamic considerations into play and
represents the Riemannian geometry induced by the second differential of Gibbs
free enthalpy $G$
\begin{equation*}
  G= \sum_{i=1}^n c_i[\ln (c_i /c_i^{\mathrm{eq}})- 1], \ \ W = \mbox{Hess} (G).
\end{equation*}
It measures the thermodynamic anisotropy of the phase space by weighting the
coordinate axis corresponding to species $i$ with the gradient $\frac{\dd
  \mu_i}{\dd c_i}=\mbox{Hess}(G)_{i,i}$ of the chemical potential
$\mu_i=\frac{\dd G}{\dd c_i}$ of species $i$ in that direction.  The
corresponding metric has been discussed in the context of an entropic scalar
product \cite{Gorban2004}. The corresponding objective function in the general
optimization problem (\ref{eq:op}) for the $W$-norm would be
\begin{equation} \label{eq:cwof} \int_0^{t_{\textrm f}} \|J_f(c) \; f \|_W \;
  \d t.
\end{equation}

Interestingly, from a different point of view the objective functionals
(\ref{eq:c2of}) and (\ref{eq:cwof}) can also be interpreted as minimizing the
length of a trajectory in a suitable Riemannian metric.  For any continuously
differentiable curve $\gamma(t)$ on a Riemannian manifold, the length $L$ of
$\gamma$ is defined as
\begin{equation}\label{eq:mincurvArclength}
  L(\gamma) = \int_{\gamma} \sqrt{g_{\gamma(t)}(\dot{\gamma}(t),\dot{\gamma}(t))} \; {\rm d}t  
\end{equation}
with $g_{\gamma(t)}$ being a scalar product defined on the tangent space of
the curve in each point. The Riemannian metric $g_{\gamma(t)}$ might be chosen
as
\begin{equation}\label{eq:mincurvMetric}
  g_{\gamma(t)}(f,f) \coloneqq f^{\textrm{T}}\;{J^{\textrm{T}}_f(c) \cdot A \cdot  J_f(c)}\; f \; = \; \| J_f(c)\;f \|_A^2,
\end{equation}
for a positive definite matrix $A$. The ``length-minimizing'' objective
functional equivalent to criterion (\ref{eq:op:of}) is now
\begin{equation}\label{eq:minMetric}
  \min\;L(\gamma)
\end{equation}
subject to constraints (\ref{eq:op:dyn})--(\ref{eq:op:pv}).  With the solution
trajectory of this problem, the ``minimum distance from equilibrium in a
kinetic sense'' can be formulated in an explicit mathematical form based on
concepts from differential geometry.  In \cite{Reinhardt2008}, a heuristic
choice for a matrix $A$ was made based on the entropy production
rate. However, the results achieved using the norm proposed in
(\ref{eq:crit_cw}) yield more accurate results for the computation of slow
attracting manifolds in chemical kinetics. Another heuristic interpretation of
(\ref{eq:c2of}) is possible based on the fact that the time-integral over a
rate of change of velocity is time-averaged velocity whose minimum is
intuitively related to the notion of a slow manifold.

\subsubsection{Geometric curvature of curves in space}\label{sss:geo_curv}
In the face of our central aim to relate the principle ``force equals
curvature'' to model reduction for chemical kinetics we refer to some
mathematical definitions and properties of curvature in this section. Various
concepts and corresponding definitions of curvature of manifolds and curves in
$\mathbb{R}^n$ can bee found in general literature on differential geometry as
e.g.\ \cite{Bloch1997, Carmo1976, Kuehnel2006}.

Let $\bar{c}:J=(0,L) \rightarrow \mathbb{R}^n$ be a curve defined on an open
interval $J \subset \mathbb{R}$ and parameterized by arc length $s$, meaning
 $ \left\|\frac{\d}{\d s}\bar{c}(s)\right\|_2=1$
($\|\cdot\|_2$ denoting the Euclidian norm). The value of
 $ \left\|\frac{\d^2}{\d s^2}\bar{c}(s)\right\|_2$
is a measure for the rate, how rapidly the curve pulls away from its tangent
in a neighborhood of $\bar{c}(s)$. That leads directly to the following
definition.

\begin{dfn}[Curvature]
  Let $\bar{c}:J \rightarrow \mathbb{R}^n$ be a curve parametrized by arc
  length $s \in J$. The number
  \begin{equation*}
    \kappa(s) \coloneqq \left\|\frac{\d^2}{\d s^2}\bar{c}(s)\right\|_2
  \end{equation*}
  is called \emph{curvature of a curve} $\bar{c}$ at $s$ (or at $\bar{c}(s)$).
\end{dfn}

However, in general trajectories in chemical composition space regarded as
curves in vector space are not parameterized in arc-length, but e.g.\ time
$t$. We want to compute the curvature for the case of an arbitrary
parametrization $t$. Let $c: I \rightarrow \mathbb{R}^n$ be a regular curve,
$I,J \subset \mathbb{R}$ open, $\varphi: J \rightarrow I$ the diffeomorphism
resulting in $\bar{c} \coloneqq c \circ \varphi$ being parametrized in arc length
with w.l.o.g.\ $\varphi(s)>0\ \forall s \in J$. Then
\begin{equation} \label{eq:cbar1diff} \frac{\d}{\d s}\bar{c}(s) = \frac{\d}{\d
    \varphi(s)}c(\varphi(s)) \frac{\d}{\d s}\varphi(s)
\end{equation}
and
\begin{equation} \label{eq:cbar2diff} \frac{\d^2}{\d s^2}\bar{c}(s) =
  \frac{\d^2}{\d \varphi(s)^2}c(\varphi(s)) \left(\frac{\d}{\d
      s}\varphi(s)\right)^2 + \frac{\d}{\d \varphi(s)}c(\varphi(s))
  \frac{\d^2}{\d s^2}\varphi(s)
\end{equation}
hold. As $\bar{c}$ is parametrized in arc length (\ref{eq:cbar1diff}) leads to
\begin{equation*}
  \frac{\d}{\d s}\varphi(s) = \frac{1}{\left\|\frac{\d}{\d \varphi(s)}c(\varphi(s))\right\|_2}.
\end{equation*}
For the second derivative of $\varphi$, that appears in (\ref{eq:cbar2diff}),
the application of the chain rule yields (with $\langle\cdot,\cdot\rangle_2$
being the Euclidian scalar product)
\begin{equation*}
  \frac{\d^2}{\d s^2}\varphi(s) =
  - \frac{\left\langle \frac{\d}{\d \varphi(s)}c(\varphi(s)), \frac{\d^2}{\d \varphi(s)^2}c (\varphi(s))\right\rangle_2} {\left\|\frac{\d}{\d \varphi(s)}c(\varphi(s))\right\|_2^4}.
\end{equation*}
Bringing the last two formulae together with (\ref{eq:cbar2diff}) and
$t=\varphi(s)$ we arrive at the formula for the curvature of $c(t)$:
\begin{equation} \label{eq:local_curv} \kappa(t) \coloneqq \Phi_{\textrm{C}}(c(t)) \coloneqq
  \left\| \frac{\ddot{c}(t)}{\|\dot{c}(t)\|_2^2} - \langle \dot{c}(t),
    \ddot{c}(t)\rangle_2 \frac{\dot{c}(t)}{\|\dot{c}(t)\|_2^4} \right\|_2.
\end{equation}

Recalling the discussions in the last section, an alternative optimization
criterion for (\ref{eq:op:of}) could be the curvature
(\ref{eq:local_curv}). In this context the total (integrated) curvature should
be the objective function in (\ref{eq:op}):
\begin{equation*}
  \kappa_{\rm tot} \coloneqq \int_{l(0)}^{l(t_{\text{f}})} \kappa(s) \d s,
\end{equation*}
which can be expressed in time-parameterization as
\begin{equation} \label{eq:curv_tot}
    \kappa_{\rm tot}
    = \int_{0}^{t_{\mathrm{f}}} \kappa(t) \| \dot{c}(t) \|_2 \; \d t
    = \int_{0}^{t_{\mathrm{f}}} \left\| \frac{\ddot{c}(t)}{\|\dot{c}(t)\|_2}
      - \langle \dot{c}(t), \ddot{c}(t)\rangle_2
      \frac{\dot{c}(t)}{\|\dot{c}(t)\|_2^3} \right\|_2 \d t
\end{equation}
with $\kappa(t)$ as in equation (\ref{eq:local_curv}).

Intuitively, the local curvature of a trajectory on its way to equilibrium in
phase space should have a peak each time it relaxes onto a lower dimensional
manifold. Therefore, also this criterion is related to the relaxation of
chemical forces in some sense.

\subsubsection{Evaluation of the Objective Functional}\label{sss:comp_of}
From a practical perspective, the computation of the Jacobian for the
expression of the different criteria is not necessary, as
 $ \ddot{c}(t) = J_f(c(t)) \; f(c(t))$
simply is a directional derivative of the ODE vector field with respect to its
own direction. This directional derivative could be evaluated using classical difference
quotients \cite{Stoer2002}, but a more appealing alternative is found in
\cite{Squire1998}. Instead of using the central difference formula
\begin{equation}
  \label{eq:centraldifference}
  F'(x_0) \approx \frac{F(x_0+\delta)-F(x_0-\delta)}{2\delta}
\end{equation}
for the approximation of the derivative of the real valued function $F(x)$,
Squire and Trapp \cite{Squire1998} suggest replacing $\delta$ with
$\mathrm{i}\delta$ ($\mathrm{i}=\sqrt{-1}$). If $F$ is an analytic function,
(\ref{eq:centraldifference}) then reads
\begin{equation}
  \label{eq:i_trick}
  F' \approx \frac{\Im[{F(x_0+\mathrm{i}\delta)}]}{\delta},
\end{equation}
with $\Im(z)$ being the imaginary part of $z$. This is called \emph{complex-step
derivative approximation}. This result is especially appealing, as (\ref{eq:i_trick}) does not contain a
subtraction and hence eliminates cancellation errors. Therefore $\delta$ can
be chosen very small, hence making higher-order terms in the Taylor expansion
negligible. For the directional derivative $\ddot{c}$ at a point $c_0$ with $c$ from
$\dot{c}=f(c)$, (\ref{eq:i_trick}) reads
\begin{equation}
  \label{eq:i_trickdirder}
  \ddot{c}\vert_{c_0} \approx \frac{\Im[f(c+\mathrm{i}\delta f(c))]}{\delta}.
\end{equation}

Compared to the use of the full Jacobian, the complexity for the evaluation of
$\ddot{c}$ can be reduced from $\mathcal{O}(n^2)$ to $\mathcal{O}(n)$ using
this complex variable approach. At the same time a high accuracy is guaranteed
by the possibility of using an extremely small $\delta$.

However, a numerical difficulty occurs for the evaluation of the objectives
(\ref{eq:cwof}) and (\ref{eq:curv_tot}). In case of (\ref{eq:cwof}) the
weights for the $W$-norm are obtained as inverted species
concentrations. Especially for radical species the denominator becomes
generally very small near chemical equilibrium resulting in numerical
instabilities. The case of (\ref{eq:curv_tot}) is even more difficult, as
negative exponents $> 1$ occur for the norm of the reaction rates. Near the
equilibrium point the reaction rates become infinitesimally small and this
results in severe numerical problems. A remedy for this problem is an
additional equality constraint. Instead of fixing the final time
$t_{\textrm{f}}$ at a large value, it can be left free in the optimization
determined by an end point constraint
\begin{equation}\label{eq:epc}
  \|f(c(t_{\textrm{f}}))\|_2 = \epsilon
\end{equation}
with a sufficiently large $\epsilon$ keeping the end point of the trajectory
away from equilibrium.

\section{Results}\label{s:res}
In this section, results for our model reduction method based on trajectory
optimization are presented. We choose three different chemical reaction
mechanisms to demonstrate its application: the Davis--Skodje model system, a
simplified reaction mechanism for the combustion of H$_2$, and a realistic
temperature dependent mechanism for ozone decomposition. For these three
mechanisms all previously discussed objective functionals (\ref{eq:c2of}),
(\ref{eq:cwof}), and (\ref{eq:curv_tot}) in the general problem (\ref{eq:op})
are tested for the purpose of numerically computing approximations of slow
attracting manifolds. The results are compared and discussed.

\subsection{The Davis--Skodje Test Problem} \label{ss:DS} The well-known
Davis--Skodje mechanism is our first test case \cite{Davis1999, Singh2002}.
\begin{equation*}
    \frac{\d y_1}{\d t}  = - y_1, \quad
    \frac{\d y_2}{\d t}  = - \gamma y_2 + \frac{(\gamma - 1)y_1 + \gamma
      y_1^2}{(1+y_1)^2},
\end{equation*}
where $\gamma > 1$ is a measure for the spectral gap or stiffness of the
system.  Typically model reduction algorithms show a good performance for
large values of $\gamma$, which represent a large gap between the time scales
of fast and slow modes. Small values of $\gamma$ impose a significantly harder
challenge on the computation of slow attracting manifolds. For reasons of
adjustable time scale separation and analytically computable slow invariant
manifold and ILDM, the Davis--Skodje model is widely used for testing
numerical model reduction approaches.

\subsubsection{Results for Different Optimization Criteria}
\begin{figure}[htbp]
  \begin{center}
    \includegraphics[width=10.cm]{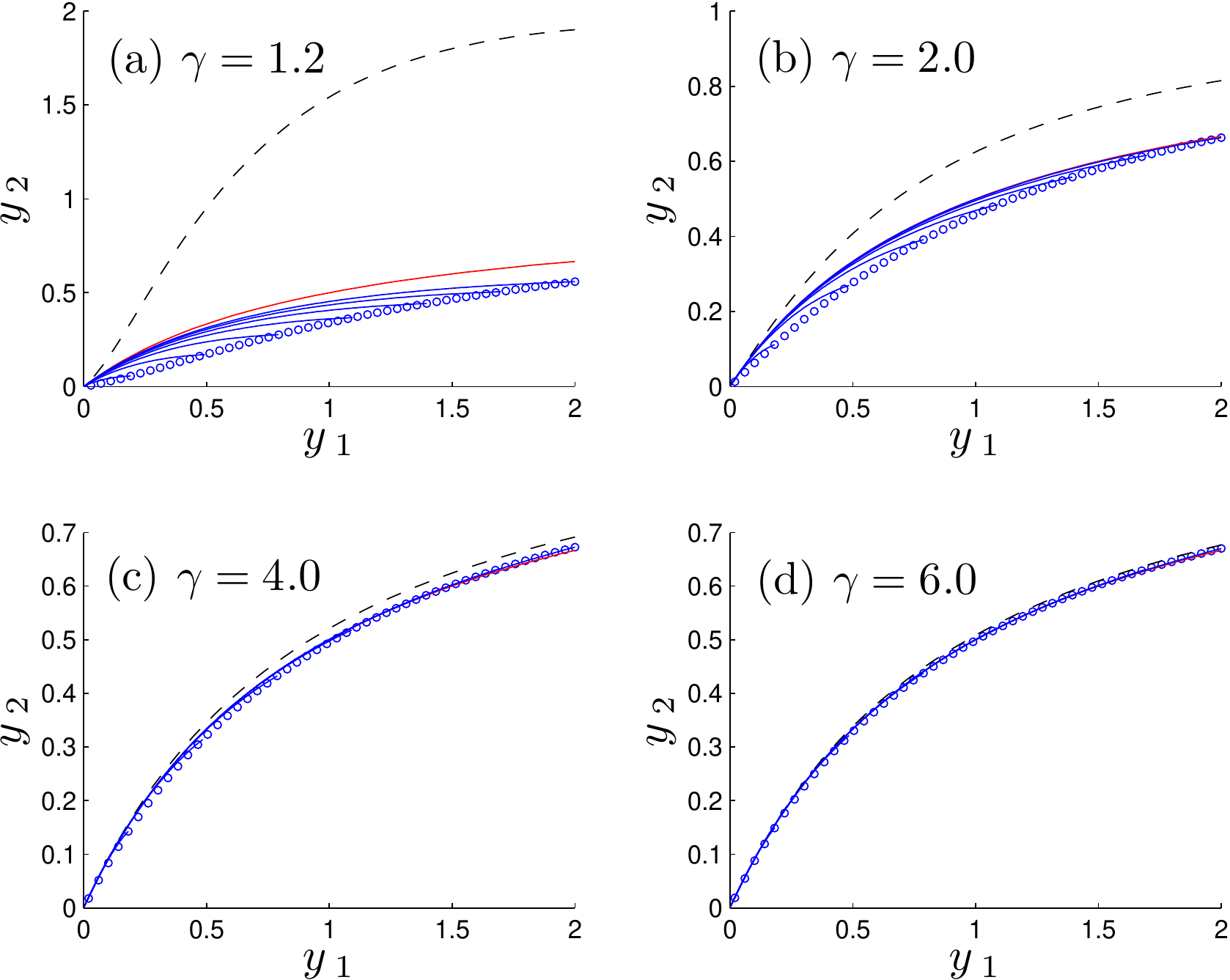}
  \end{center}
  \caption{\label{f:ds_jf2} Results for the Davis--Skodje problem with
    (\ref{eq:c2of}) as objective functional. Results for different values of
    $\gamma$ are shown. The red curve is the analytically computed SIM. The
    black dashed curve represents the analytic Maas--Pope-ILDM. The blue
    curves are trajectories numerically integrated from those initial points
    that are solutions of our optimization problem.}
\end{figure}

\begin{figure}[htbp]
  \begin{center}
    \includegraphics[width=10.cm]{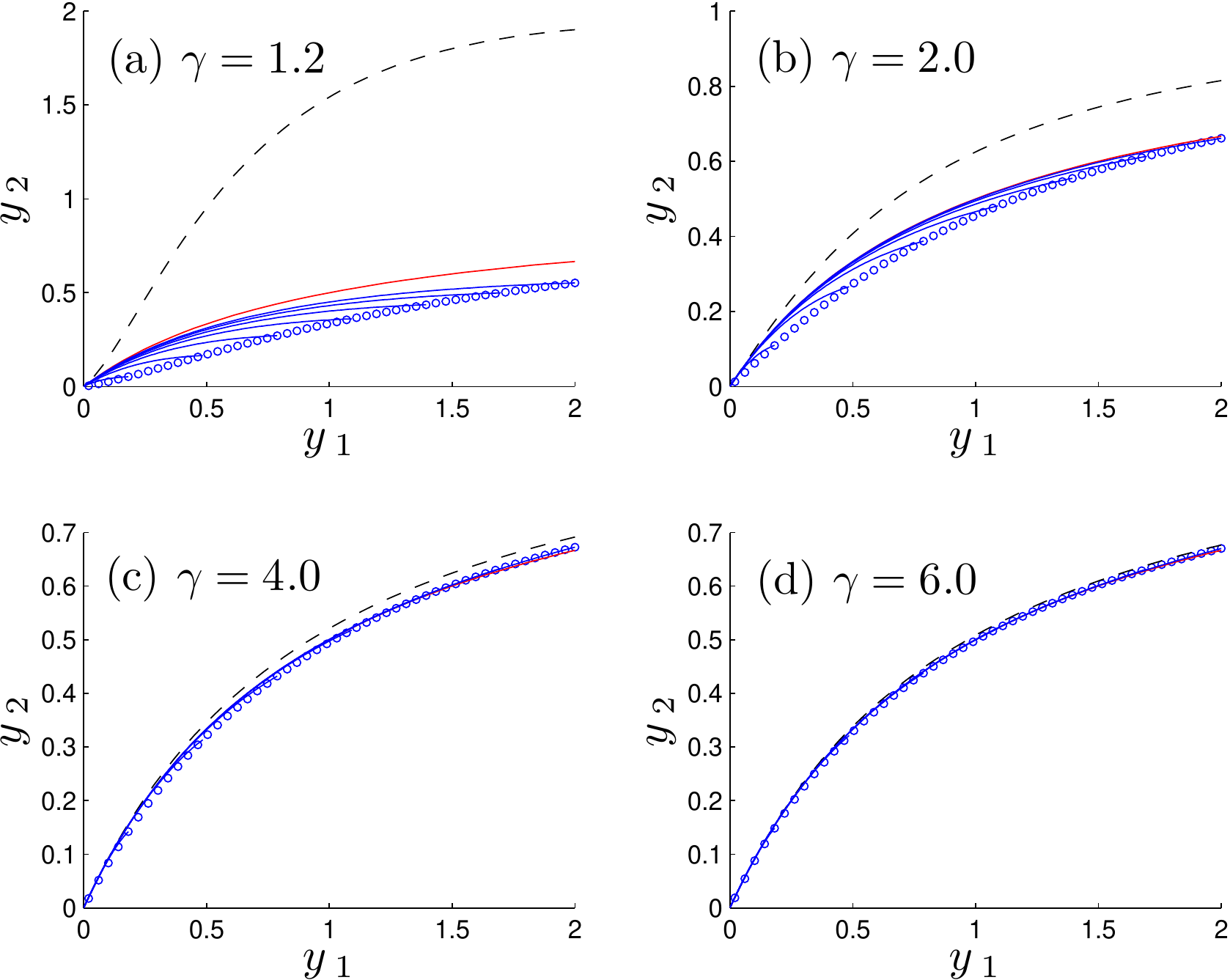}
  \end{center}
  \caption{\label{f:ds_jfw} Results for the Davis--Skodje problem with
    (\ref{eq:cwof}) as objective functional. Again, results for different
    values of $\gamma$ are shown together with the SIM (red) and the
    Maas--Pope-ILDM (black, dashed).}
\end{figure}

\begin{figure}[htbp]
  \begin{center}
    \includegraphics[width=10.cm]{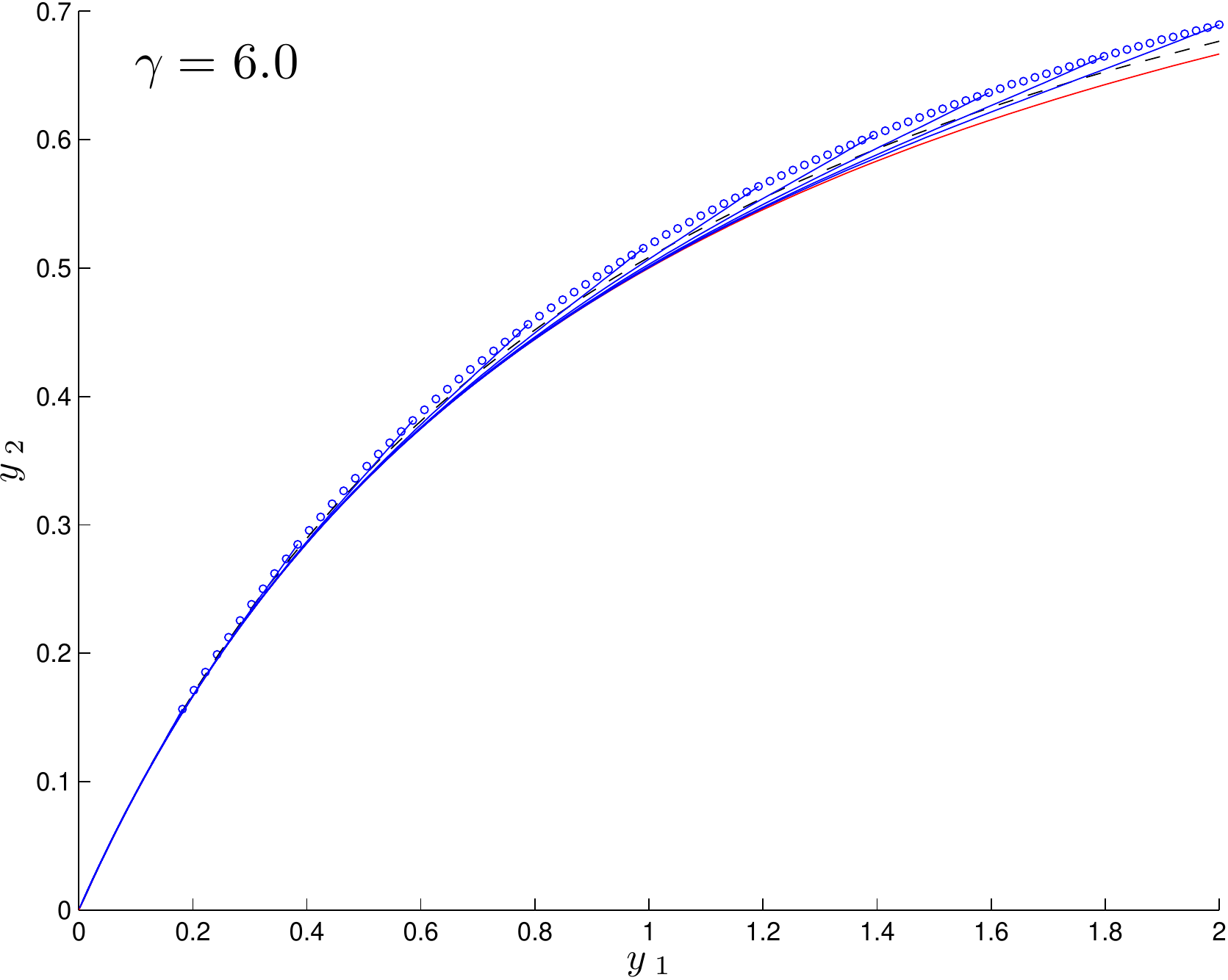}
  \end{center}
  \caption{\label{f:ds_mep} Results for the Davis--Skodje problem with
    (\ref{eq:curv_tot}) as objective function. Only results for $\gamma = 6.0$
    are shown, cf.\ the discussion in Section \ref{sss:ds_ol} including Figures
    \ref{f:ds6MEP} and \ref{f:ds2MEP}.}
\end{figure}

In Figures \ref{f:ds_jf2}, \ref{f:ds_jfw}, and \ref{f:ds_mep} results for
criteria (\ref{eq:c2of}), (\ref{eq:cwof}), and (\ref{eq:curv_tot})
respectively are depicted.

We will refer to criterion (\ref{eq:crit_c2}) as criterion A in the
following. The dependence of the accuracy of the computed slow attracting
manifold on the stiffness parameter $\gamma$ becomes obvious. In all cases,
the value of $y_1$ is fixed as reaction progress variable. The optimization
problem is solved repeatedly for different values for $y_1$. For large and
moderate values of the stiffness parameter, good approximations of the SIM
(red) are achieved. For $\gamma = 1.2$ the results become more inaccurate. For
comparison the Maas--Pope-ILDM is plotted as dashed black line, it can be
computed analytically for the Davis--Skodje model \cite{Davis1999}. The
computational effort for one solution of the optimization problem is on
average about twelve SQP-iterations (cf.\ Section \ref{ss:numerics}) which
takes about ten seconds in total on a single core Intel$^\circledR$
Pentium$^\circledR$ 4 (3 GHz)-machine with 2 GB memory. Of course, the
convergence time (not the result) depends on the initial values chosen to
start the numerical algorithm. We use a non-equidistant multiple shooting grid
with twenty intervals. The number of evaluations of the function $f$ is of
order $10^5$, the order of the number of matrix factorizations is $10^4$. Due
to the parametric optimization strategy and application of initial value
embedding \cite{Diehl2005, Diehl2002b} pointed out in Section
\ref{ss:numerics}, every follow-up solution of a neighboring problem with
slightly different values for the reaction progress variables needs only three
to five SQP-iterations.

For the second criterion (\ref{eq:cwof}) -- denoted B in the following -- the
results look very similar (see Fig.~\ref{f:ds_jfw}). As in the Davis--Skodje
model the values for the variables $y_1$ and $y_2$ are of the same order, this
is obvious considering the scaling in criterion (\ref{eq:crit_cw}). Results
for criterion C (\ref{eq:curv_tot}), the total curvature, are shown in
Figure \ref{f:ds_mep}. With this criterion, we could not obtain reasonable
results for values of $\gamma < 6.0$ which will be further explained in the
next section.

\subsubsection{Optimization Landscapes} \label{sss:ds_ol} Since the
Davis--Skodje test problem consists of only two variables, the structure of
optimization landscapes can easily be visualized. To compute these
landscapes, the initial values of both variables are varied over a fixed
range and for the trajectories starting in each of these pairs of initial
values, the values of (\ref{eq:c2of}), (\ref{eq:cwof}), or
(\ref{eq:curv_tot}) respectively are calculated. These objective functional
values are depicted as a function of the initial values of the corresponding
trajectories. Calculations are performed and plots are generated using
MATLAB$^\circledR$.

\begin{figure}[htbp]
 \begin{center}
  \subfigure[ Results with $\gamma = 6.0$ corresponding to results depicted in
    Fig.~\ref{f:ds_mep}.]{\includegraphics[width=0.4\textwidth]{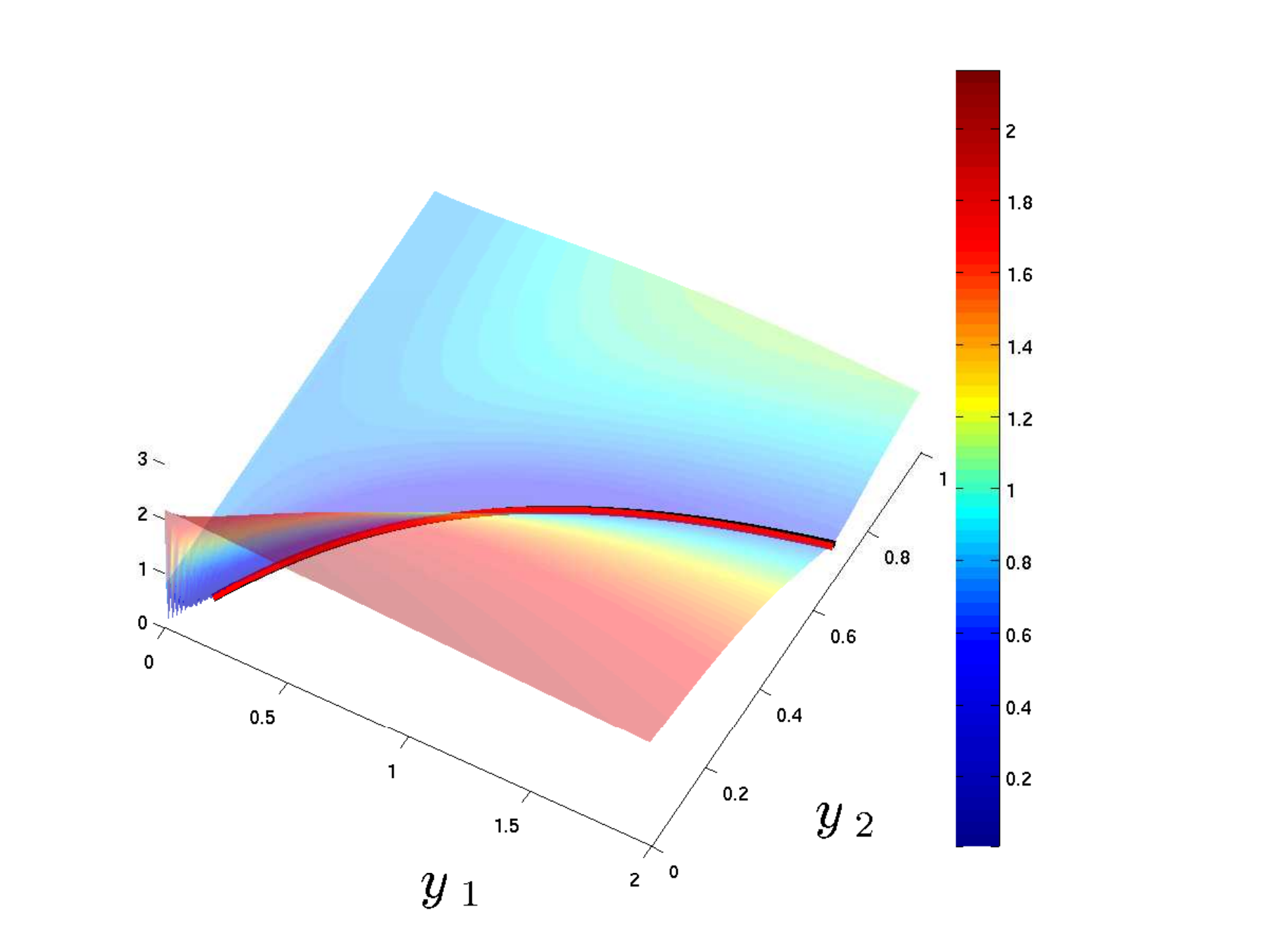}\label{f:ds6MEP}} \hspace*{0.1\textwidth}
  \subfigure[ Results with $\gamma = 2.0$. No minimum near the SIM can be
    found.]{\includegraphics[width=0.4\textwidth]{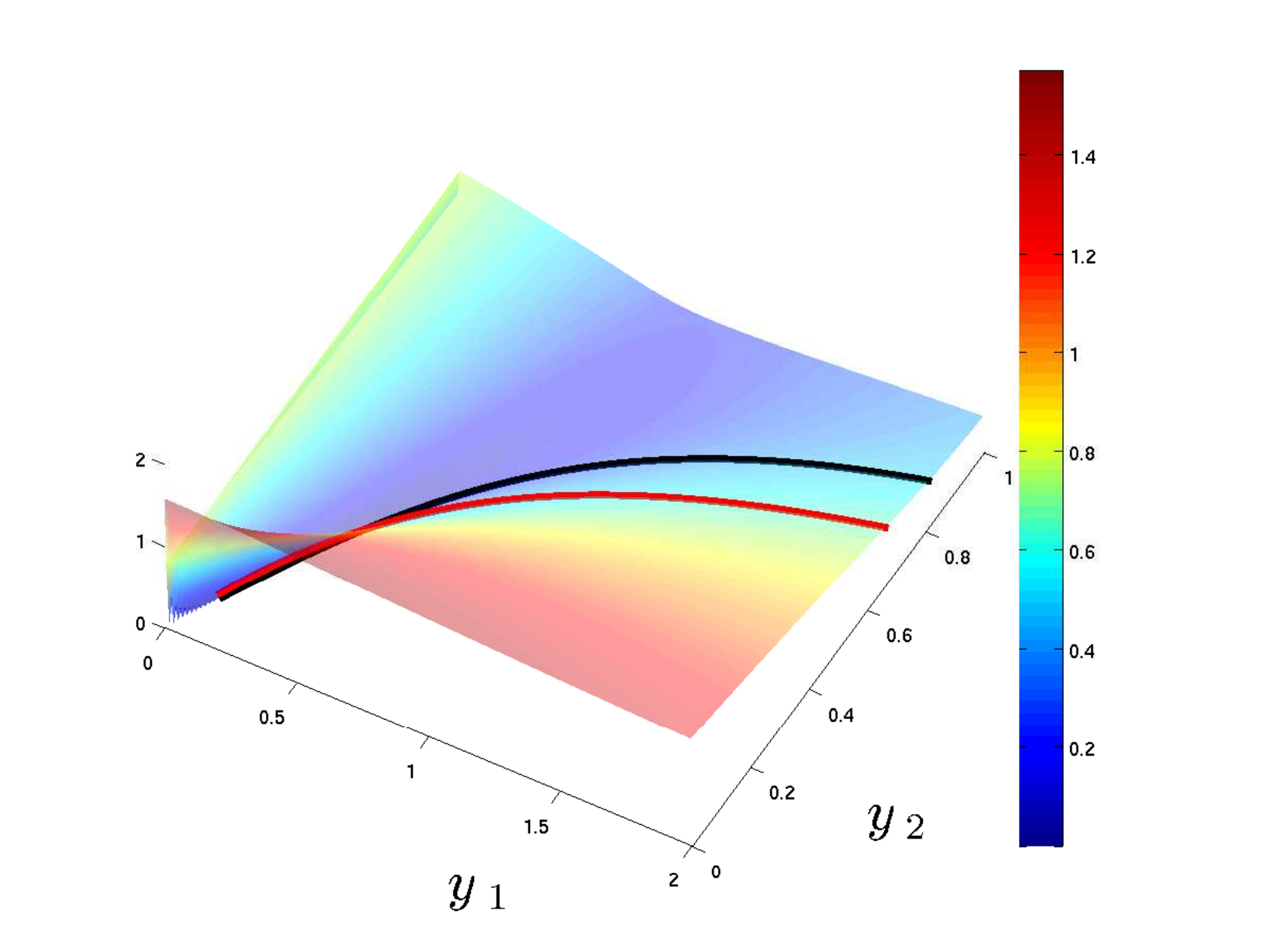}\label{f:ds2MEP}}
 \end{center}
 \caption{ Optimization landscape for the
    Davis--Skodje-model. The integrated (total) curvature
    (\ref{eq:curv_tot}) is plotted on the $z$-axis and coded in color for
    illustration reasons. The analytically computed SIM (red) and the
    analytically given Maas--Pope-ILDM (black) are projected onto the
    landscape.}
\end{figure}

For the Davis--Skodje problem, we restrict ourselves to the illustration of
criterion C, where we did not achieve satisfying results for small values of
$\gamma$ (see previous subsection). In Figure \ref{f:ds6MEP} the results for
$\gamma = 6.0$ are shown. Additionally the SIM (red) and the ILDM (black)
are projected onto the optimization landscape. Remember, that in our
approach for a fixed value of $y_1$ the minimum of the objective value
identifies the corresponding initial value of $y_2$.

An optimization landscape for $\gamma = 2.0$ and criterion C is shown in
Figure \ref{f:ds2MEP}. In this case, the criterion fails. No minimum can be
found in the neighborhood of the SIM. The reason for this is that obviously
for small time scale separation the relation between geometric and kinetic
properties of the trajectories pointed out in Section~\ref{sss:geo_curv}
becomes weaker and linear segments of trajectories are preferred in the
optimization objective functional against a slow attracting manifold (see
``valley'' in Fig.~\ref{f:ds2MEP}). As explained in
Section~\ref{sss:geo_curv}, criterion C measures the curvature generated
through relaxation of a trajectory onto a slow attracting manifold. Thus, if
relaxation is weak, the criterion becomes inaccurate.

\subsection{Model Hydrogen Combustion Reaction Mechanism}
In this section we consider a small test mechanism, which has been used for
model reduction purposes in \cite{Chiavazzo2008,Gorban2004,Reinhardt2008}. It
consists of six chemical species involved in six (in each case forward and
backward) elementary reactions involving two element mass conservation
relations for hydrogen and oxygen (cf.\ Table \ref{mech:H2comb}).

\begin{table}[ht]
  \caption{\label{mech:H2comb} Simple hydrogen combustion test mechanism from \cite{Gorban2004}. Forward and backward rate constants are given temperature-independently.}
  \centerline{
    \begin{tabular}{@{}lclrr@{}} 
      \toprule 
Reaction               &                      &                        & $k_{+}$  & $k_{-}$   \\
\midrule
$\spec{H_2}$           & $\rightleftharpoons$ & $\spec{2}\,\spec{H}$   & $2.0$    & $216.0$   \\
$\spec{O_2}$           & $\rightleftharpoons$ & $\spec{2}\,\spec{O}$   & $1.0$    & $337.5$   \\
$\spec{H_2O}$          & $\rightleftharpoons$ & $\spec{H} + \spec{OH}$ & $1.0$    & $1400.0$  \\
$\spec{H_2}+ \spec{O}$ & $\rightleftharpoons$ & $\spec{H} + \spec{OH}$ & $1000.0$ & $10800.0$ \\
$\spec{O_2}+ \spec{H}$ & $\rightleftharpoons$ & $\spec{O} + \spec{OH}$ & $1000.0$ & $33750.0$ \\
$\spec{H_2}+ \spec{O}$ & $\rightleftharpoons$ & $\spec{H_2O}$          & $100.0$  & $0.7714$  \\
      \bottomrule
    \end{tabular}
  }
\end{table}

With the mass conservation relations
\begin{equation*}
  \begin{split}
    2  \: c_{\text{H}_2} + 2 \: c_{\text{H}_2 \text{O}} + c_{\text{H}} + c_{\text{OH}} & = C_1 \\
    2 \: c_{\text{O}_2} + c_{\text{H}_2 \text{O}} + c_{\text{O}} +
    c_{\text{OH}} & = C_2
  \end{split}
\end{equation*} 
this mechanism yields a system with four degrees of freedom. For our
computations $C_1 = 2.0$ and $C_2 = 1.0$ were chosen.

\subsubsection{One-dimensional Manifolds}
We first present results for the computation of one-dimensional manifolds in
composition space.  The value of $c_{\text{H}_2 \text{O}}$ serves as reaction
progress variable. It is varied between 0.05 and 0.65. We present and compare
results for the three different optimization criteria introduced in
Section~\ref{ss:optimcrit}. In general, for testing accuracy of model reduction approaches it is very
difficult to provide qualitative and quantitative measures of accuracy of
numerical results. The reason for this is the nonavailability of
(analytical) expressions for a SIM. Commonly we consider ``eye inspection''
of trajectory bundling behavior as well as consistency (invariance) as a
qualitative measure of accuracy.

\begin{figure}
  \begin{center}
    \includegraphics[width=10.cm]{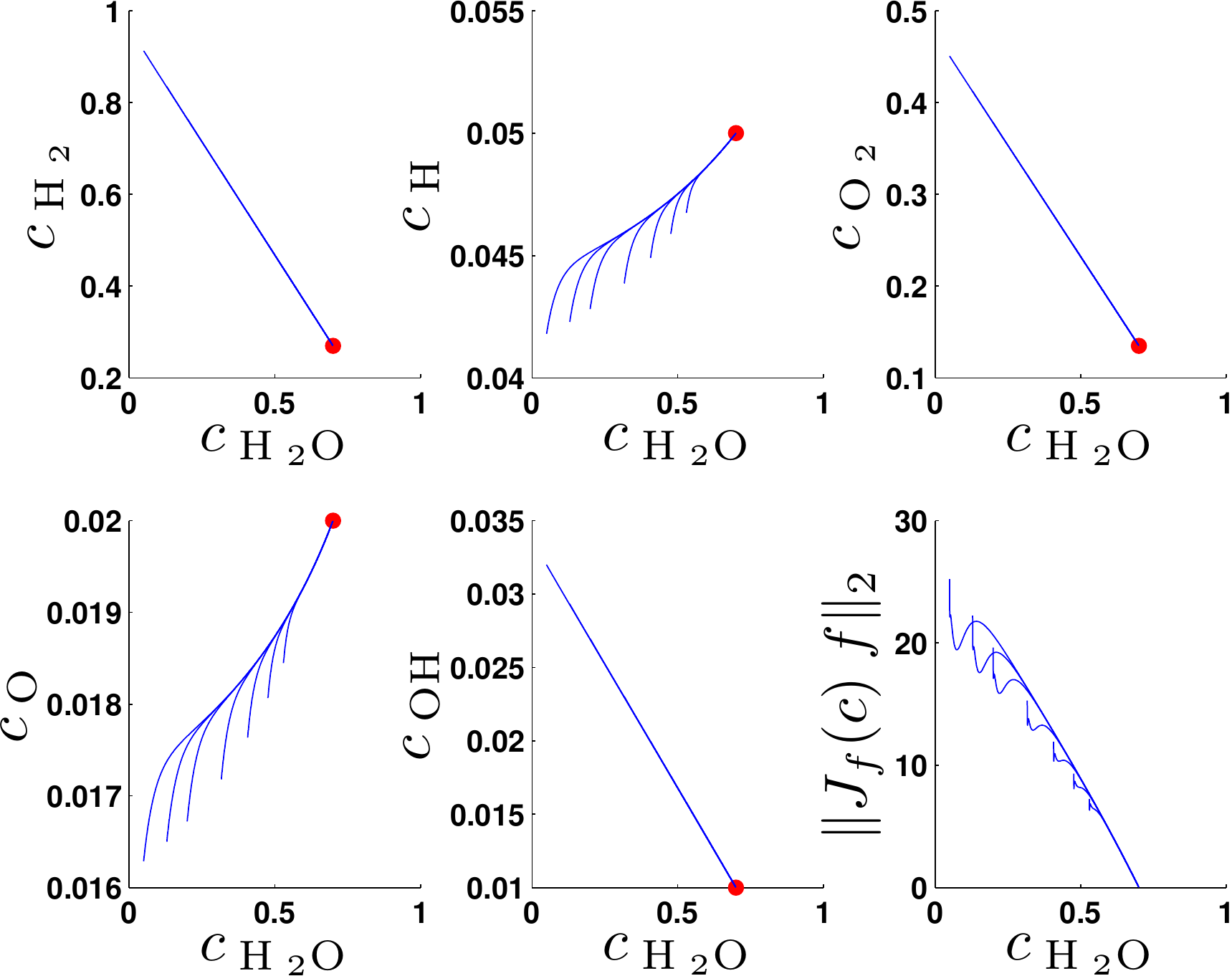}
  \end{center}
  \caption{\label{f:H2C6Jf2_1d} Results for the hydrogen combustion mechanism
    with criterion A.\@ The free variables and the integrand of the objective
    functional are plotted versus the reaction progress variable
    $c_{\mathrm{H}_2 \mathrm{O}}$. For different values of $c_{\mathrm{H}_2
      \mathrm{O}}$ the optimization problem is solved and the evolution of the
    resulting trajectories (blue curves) towards equilibrium (red dot) is
    shown. Especially for the radical species concentrations $c_{\mathrm{H}}$
    and $c_{\mathrm{O}}$ criterion A turns out to be inconsistent.}
\end{figure}

\begin{figure}
  \begin{center}
    \includegraphics[width=10.cm]{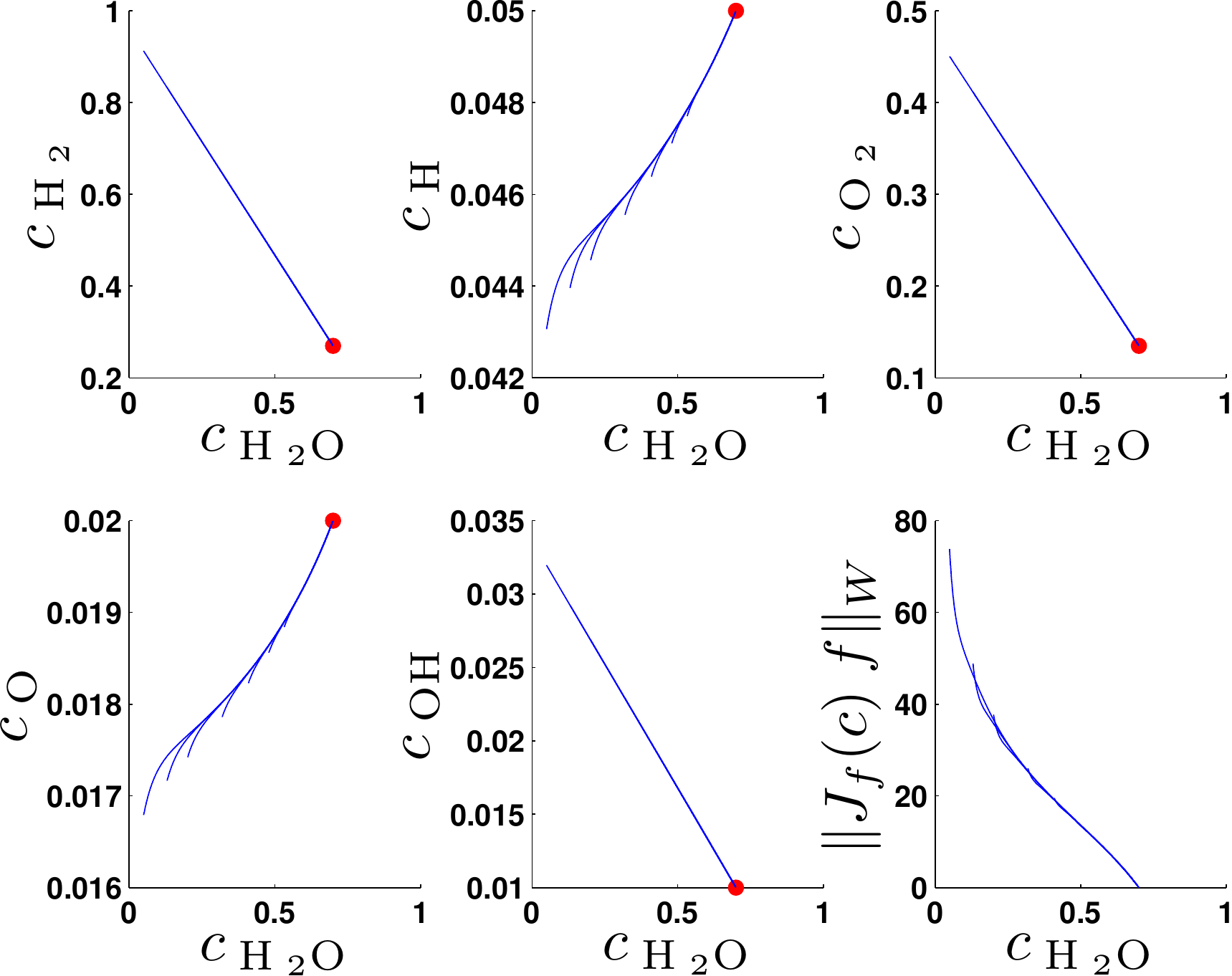}
  \end{center}
  \caption{\label{f:H2C6JfW_1d} Results for the hydrogen combustion mechanism
    with criterion B.}
\end{figure}

\begin{figure}
  \begin{center}
    \includegraphics[width=10.cm]{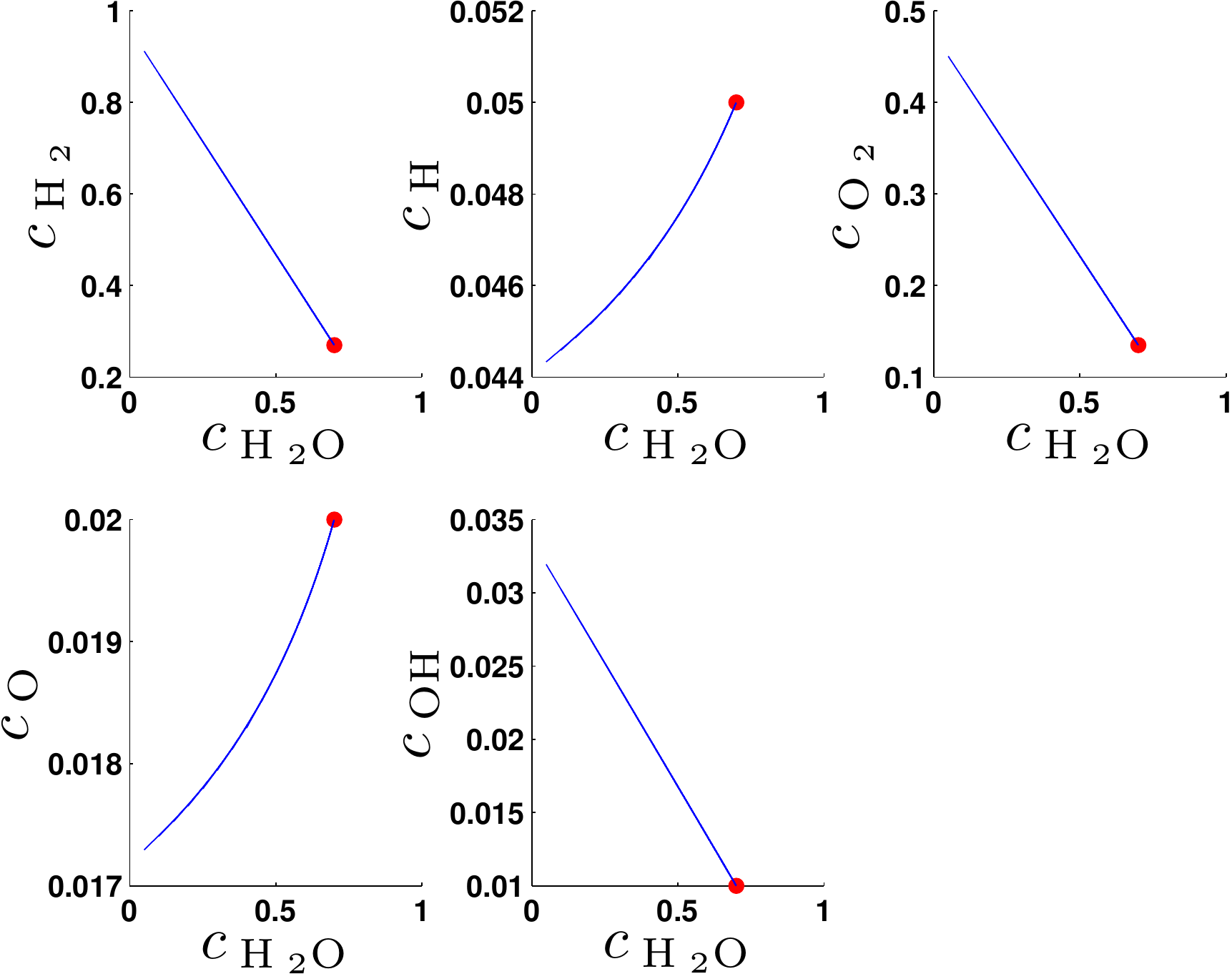}
  \end{center}
  \caption{\label{f:H2C6MEP_1d} Results for the hydrogen combustion mechanism
    with criterion C.}
\end{figure}

In Figure \ref{f:H2C6Jf2_1d}, results for criterion A are shown. The values of
the free variables computed in the optimization are plotted versus the value
of $c_{\text{H}_2 \text{O}}$. Especially for the radical species
concentrations $c_{\text{H}}$ and $c_{\text{O}}$ criterion A turns out to be
inconsistent to some extent in the sense of Definition \ref{def:consis}.
Figure \ref{f:H2C6JfW_1d} shows the results for the weighted criterion B.\@ A
significant improvement towards better consistency is achieved.
The third criterion C which failed in case of the Davis--Skodje test problem
here performs best, cf.\ Figure \ref{f:H2C6MEP_1d}. The results are nearly
consistent. For criterion C the additional equality constraint (\ref{eq:epc}),
proposed in Section~\ref{sss:comp_of} has been used to prevent numerical
instabilities near the equilibrium point. The computational effort for a 1-D
manifold is in the order of seconds.

\subsubsection{Two-dimensional Manifolds}
As the hydrogen combustion model has four degrees of freedom, also
two-dimensional manifolds can be constructed. In the presented examples,
$c_{\text{H}_2 \text{O}}$ and $c_{\text{H}_2}$ serve as reaction progress
variables. We present consistency tests plotted in two dimensions and finally
show three-dimensional plots of the computed two-dimensional manifold and the
relaxation of trajectories started from arbitrary initial values onto this 2-D
manifold.

\begin{figure}
  \begin{center}
    \includegraphics[width=10.cm]{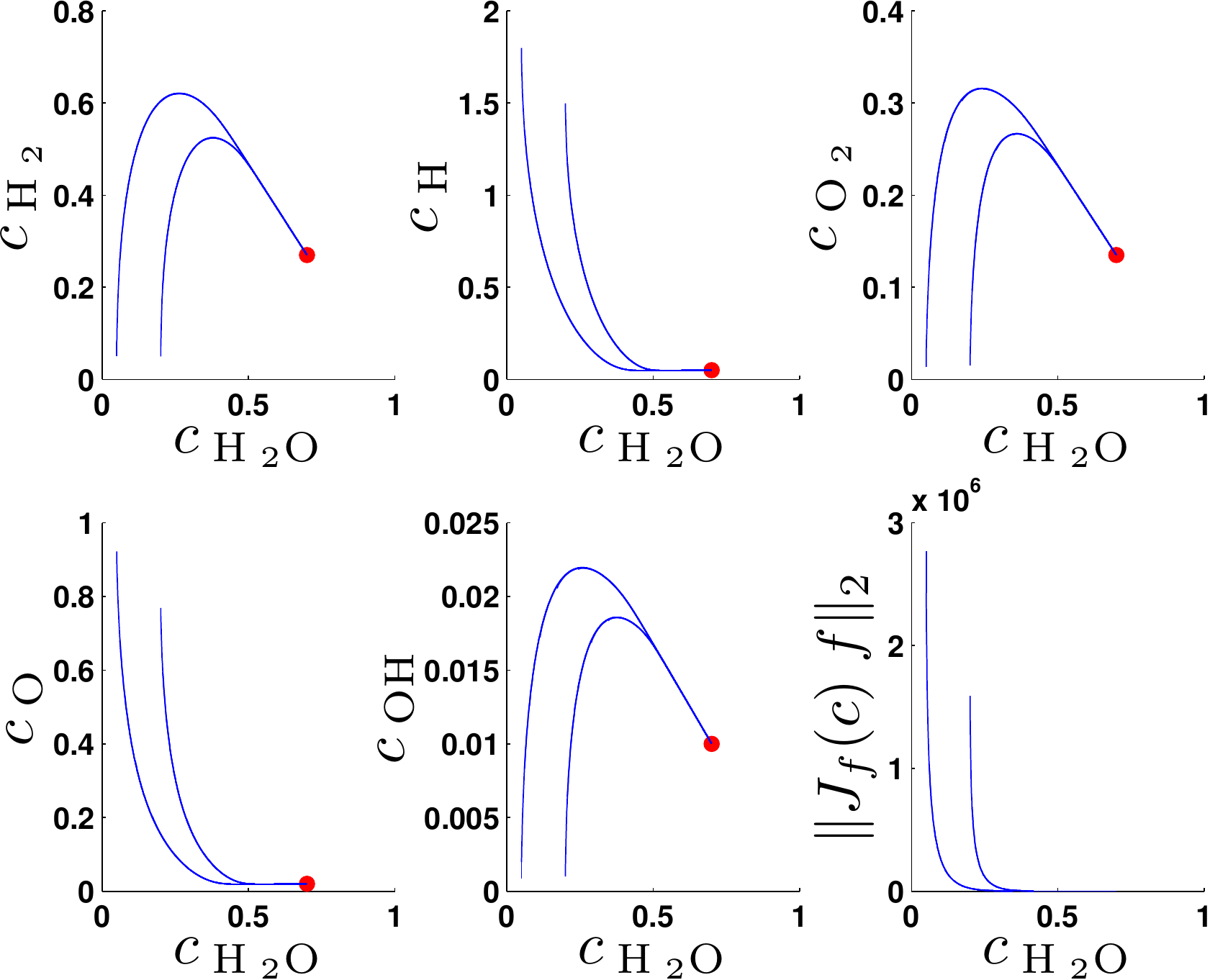}
  \end{center}
  \caption{\label{f:H2C6Jf2_2d} Results for a two-dimensional manifold for the
    hydrogen combustion mechanism with criterion A.\@ The values of the
    progress variables are fixed to $0.05$ and $0.2$ in case of
    $c_{\mathrm{H}_2 \mathrm{O}}$ and to $0.05$ in case of
    $c_{\mathrm{H}_2}$. Consistency tests are performed.}
\end{figure}

\begin{figure}
  \begin{center}
    \includegraphics[width=10.cm]{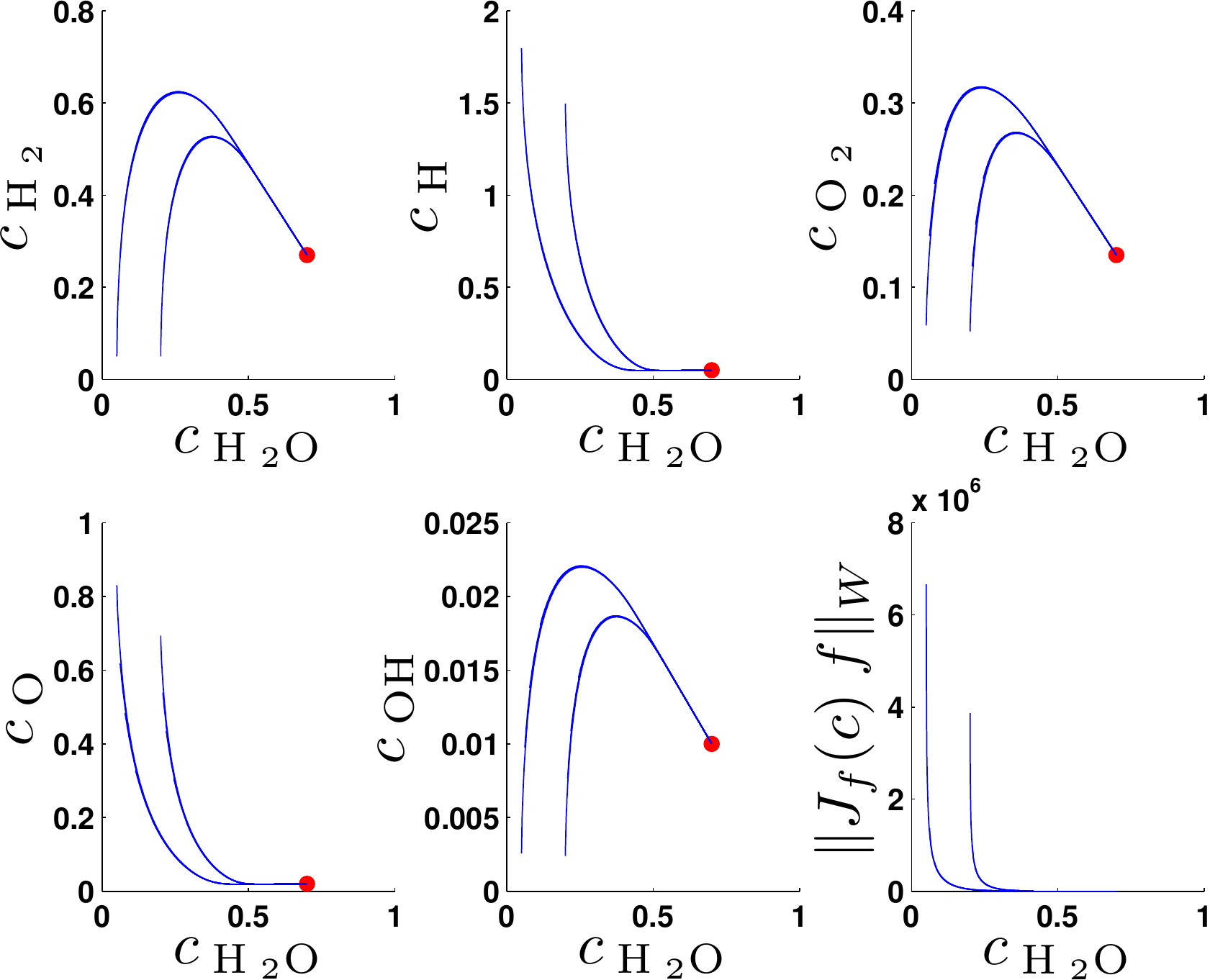}
  \end{center}
  \caption{\label{f:H2C6JfW_2d} Results for a two-dimensional manifold for the
    hydrogen combustion mechanism with criterion B.\@ The results are
    comparable to the results for criterion A, cf.\
    Figure~\ref{f:H2C6Jf2_2d}.}
\end{figure}

Figure \ref{f:H2C6Jf2_2d} refers to criterion A.\@ The resulting trajectories
are plotted. After some time $t_1$ the values of the progress variables are
fixed and the problem is solved again for testing consistency by eye
inspection which turns out to be quite accurate. For criterion B (see
Fig.~\ref{f:H2C6JfW_2d}), comparably good consistency is observed.

\begin{figure}
  \begin{center}
    \includegraphics[width=10.cm]{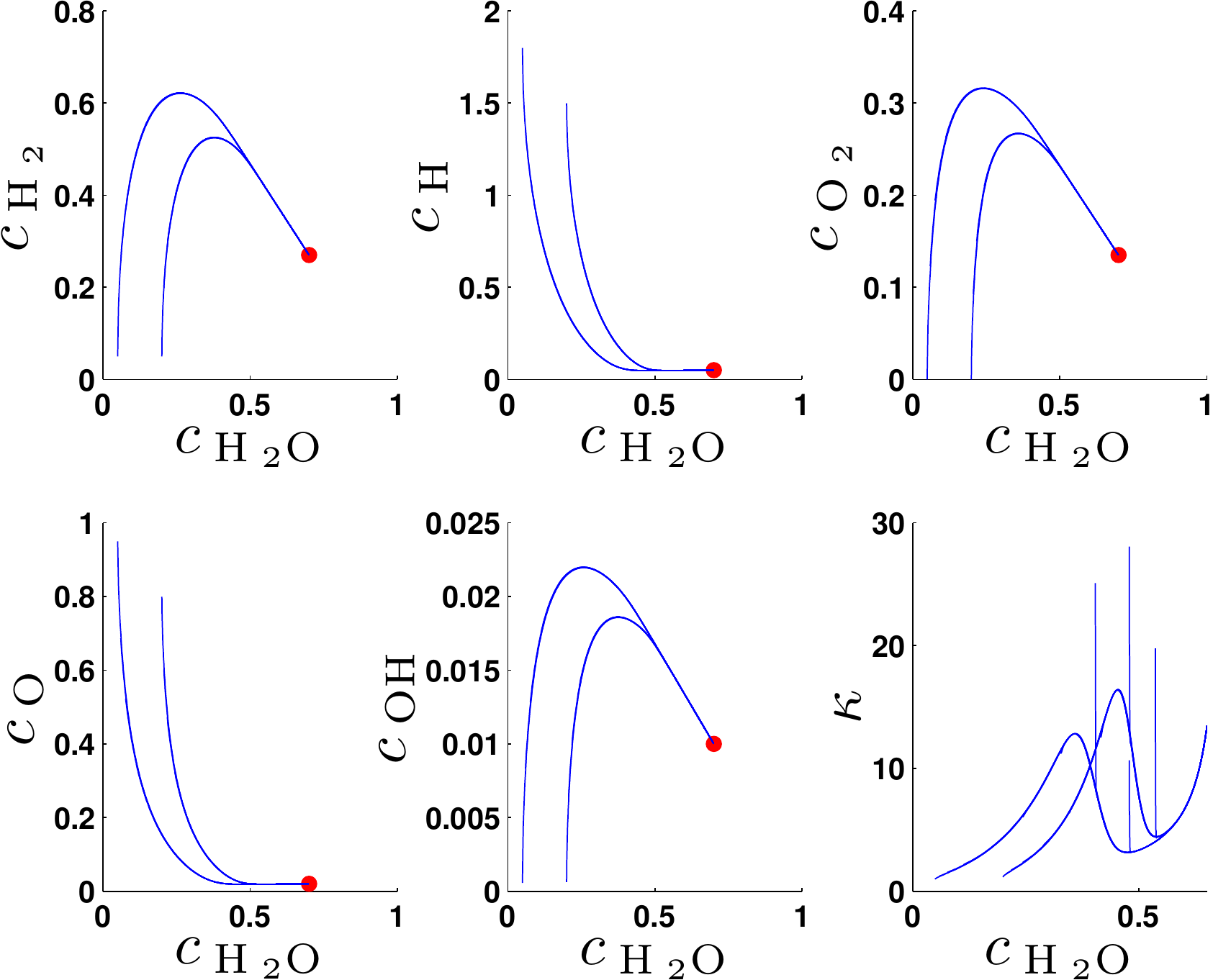}
  \end{center}
  \caption{\label{f:H2C6MEP_2d} Results of the two-dimensional manifold for
    the hydrogen combustion mechanism using criterion C.}
\end{figure}

Results for the third criterion C are depicted in Figure
\ref{f:H2C6MEP_2d}. Here the peaks in the curvature during relaxation onto a
lower-dimensional manifold mentioned in Section~\ref{sss:geo_curv} become
obvious.

\begin{figure}
  \begin{center}
    \includegraphics[width=10.cm]{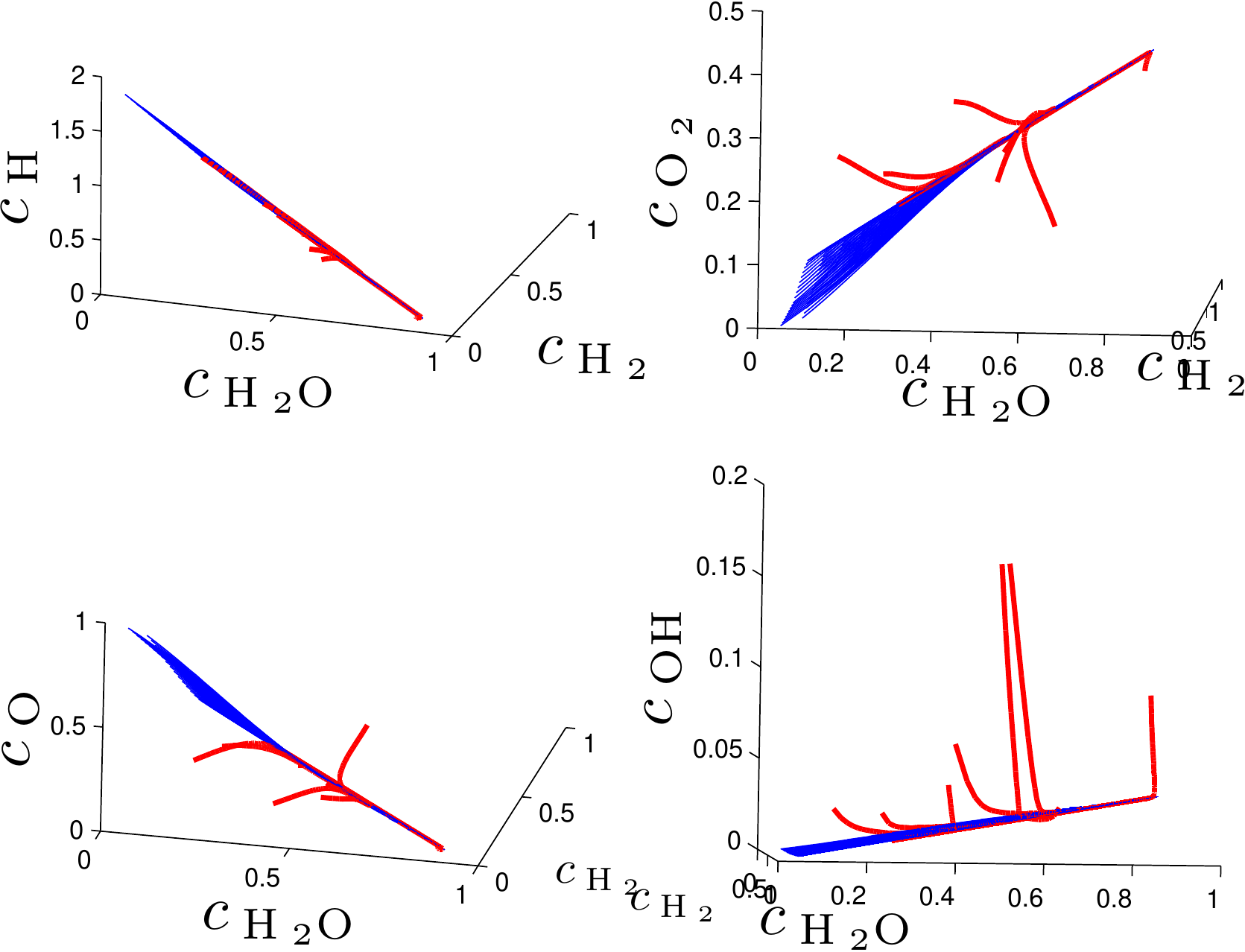}
  \end{center}
  \caption{\label{f:H2C6MEP_3d} Three-dimensional plots of the two-dimensional
    manifold for the hydrogen combustion mechanism. The free variables are
    plotted versus the progress variables. The manifold is spanned by
    trajectories (blue) computed for initial fixed values $c_k^0$ with
    criterion C as objective functional. Arbitrary trajectories (red)
    fulfilling the element mass conservation are computed to visualize their
    relaxation onto the manifold.}
\end{figure}

For visualization of the two-dimensional manifold, three-dimensional cuts of
six-dimensional composition space are plotted in Figure
\ref{f:H2C6MEP_3d}. The remaining free variables are plotted versus the
reaction progress variables. Arbitrary trajectories relax on the 2-D manifold
spanned by the computed trajectories. The computational effort for a 2-D
manifold is in the order of some minutes.

\subsection{Ozone Decomposition Reaction Mechanism}
Our last test case is a three component ozone decomposition mechanism (see
Table \ref{mech:ozone}) taken from \cite{Maas1989}. It has been chosen to
demonstrate the performance of our method taking temperature dependence via
Arrhenius kinetics into account. Many model reduction approaches explicitly
based on time scale separation fail when the spectral gap between fast and
slow modes becomes too small which is often the case for low temperatures.
\begin{table}[ht]
  \caption{\label{mech:ozone}Ozone decomposition mechanism from \cite{Maas1989}. Rate coefficient $k=AT^b\exp(-E_\mathrm{a}/RT)$. Collision efficiencies in reactions including $\mathrm{M}$: $f_{\rm O}=1.14, f_{{\rm O}_2} = 0.40, f_{{\rm O}_3} = 0.92$.}
  \centerline{
    \begin{tabular}{@{}lclrrr@{}} 
      \toprule 
      Reaction && & $A$ / $\text{cm}, \text{mol}, \text{s}$ & $b$ & $E_\mathrm{a}$ / $\frac{\text{kJ}}{\text{mol}}$ \\  \midrule
      O + O + M & $\rightarrow$ & O$_2$ + M &  $2.90\times 10^{17}$ & $-1.0$ & $0.0$  \\
      O$_2$ + M & $\rightarrow$ & O + O + M & $6.81\times 10^{18}$ & $-1.0$ & $496.0$ \\
      O$_3$ + M & $\rightarrow$ & O + O$_2$ + M & $9.50\times 10^{14}$ & $0.0$ &  $95.0$ \\
      O  + O$_2$ + M & $\rightarrow$ & O$_3$ + M & $3.32\times 10^{13}$ & $0.0$ & $-4.9$ \\
      O  + O$_3$ & $\rightarrow$ & O$_2$ + O$_2$ & $5.20\times 10^{12}$ & $0.0$ & $17.4$\\
      O$_2$ + O$_2$ & $\rightarrow$ & O + O$_3$ & $4.27\times 10^{12}$ & $0.0$ & $413.9$\\
      \bottomrule
    \end{tabular}
  }
\end{table}

The ozone decomposition mechanism involves the element mass conservation
relation
\begin{equation*}
  c_{\text{O}} + 2\:c_{\text{O}_2} + 3\:c_{\text{O}_3} = C
\end{equation*}
leaving a system with two degrees of freedom. We choose without loss of
generality $C=1$.

\subsubsection{Results for Different Optimization Criteria}\label{sss:o3or}
\begin{figure}
  \begin{center}
    \includegraphics[width=10.cm]{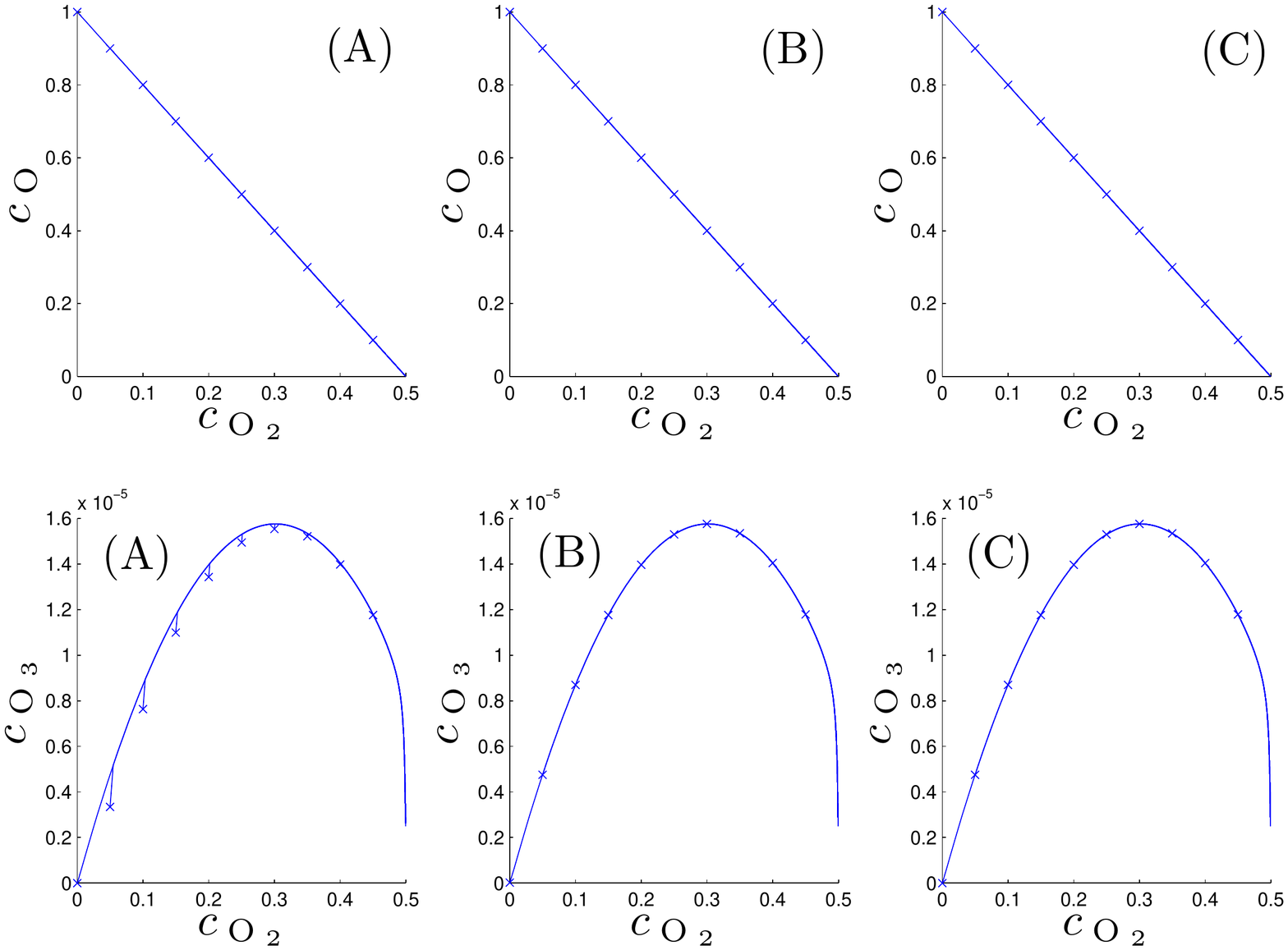}
  \end{center}
  \caption{\label{f:o3_100} Results for the ozone decomposition mechanism at a
    temperature of $T=1000~{\mathrm{K}}$ with the three different criteria A
    (left), B (middle), and C (right). The free variables are plotted versus
    $c_{\mathrm{O}_2}$ as reaction progress variable. The optimization problem
    is solved several times with different values of $c_{\mathrm{O}_2}(0)$
    (x-marks depicting the initial values of the optimal solution
    trajectories). The solution trajectories starting at the blue x-marks are
    shown on their way to equilibrium ($c_{\mathrm{O}_2}^{\mathrm{eq}} =
    0.5$). All criteria work reasonably well, but criterion A is worse
    concerning consistency.}
\end{figure}

The ozone decomposition model has two degrees of freedom and we compute
one-dimensional manifolds. The results for different temperatures are
compared. Consistency tests are performed as in the previous sections.  For
criterion A at $T=1000~{\textrm{K}}$ (Fig.~\ref{f:o3_100}), a short relaxation
phase of the computed trajectories can be observed, indicating
inconsistency. The other criteria yield more consistent results and seem to be
of comparable quality for approximation of a slow attracting manifold.

\begin{figure}
  \begin{center}
    \includegraphics[width=10.cm]{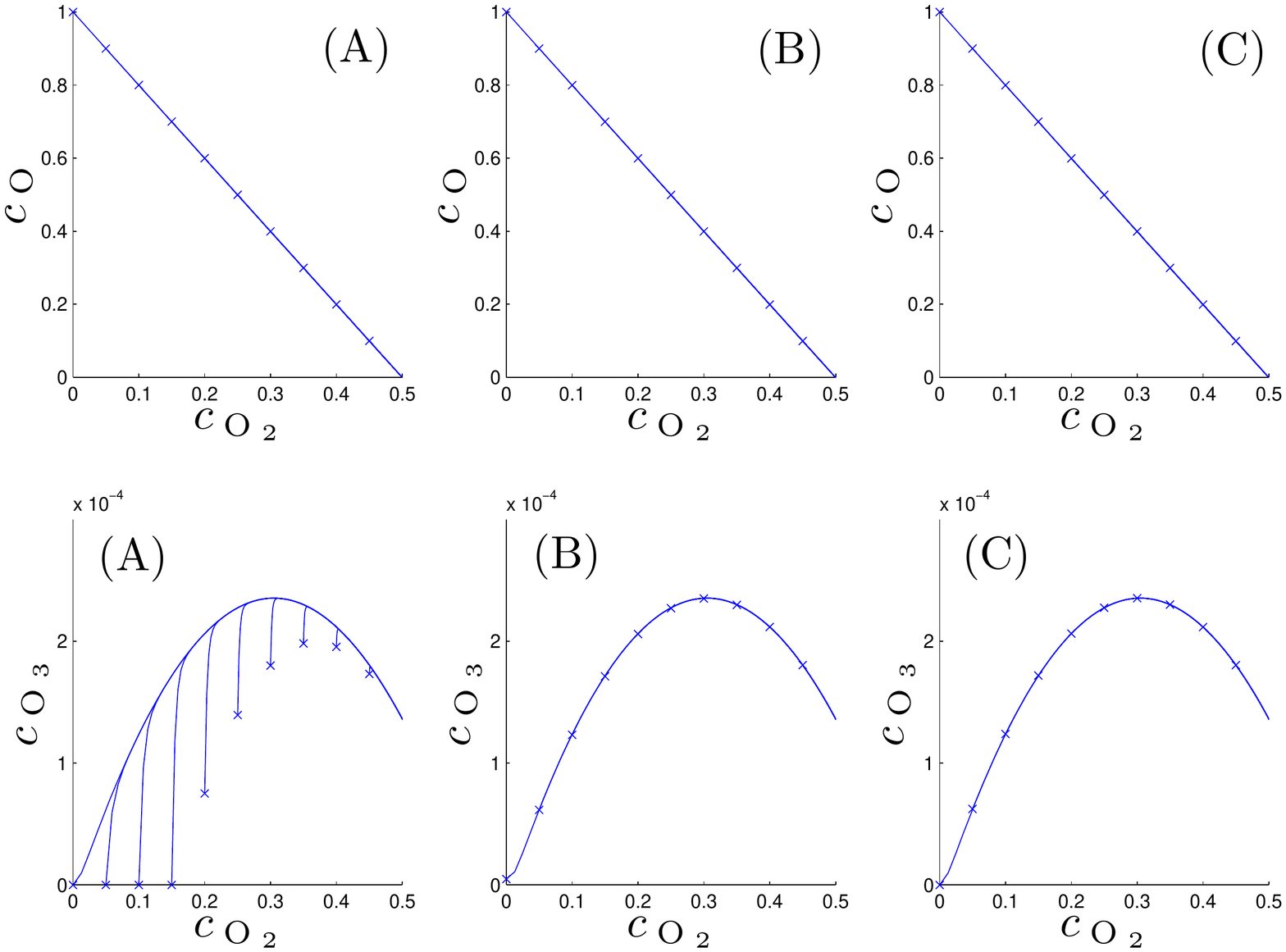}
  \end{center}
  \caption{\label{f:o3_50} Results for the ozone decomposition mechanism at a
    temperature of $T=500~{\mathrm{K}}$ with the three different criteria
    arranged as in Figure \ref{f:o3_100}. Criteria B and C perform well,
    whereas criterion A obviously fails.}
\end{figure}

For a lower temperature of $T=500~{\textrm{K}}$ (Figure \ref{f:o3_50}),
criterion A obviously fails for the ``relatively small'' absolute values of
$c_{\text{O}_3}$, whereas criteria B and C yield good approximations of the
SIM.

\begin{figure}
  \begin{center}
    \includegraphics[width=10.cm]{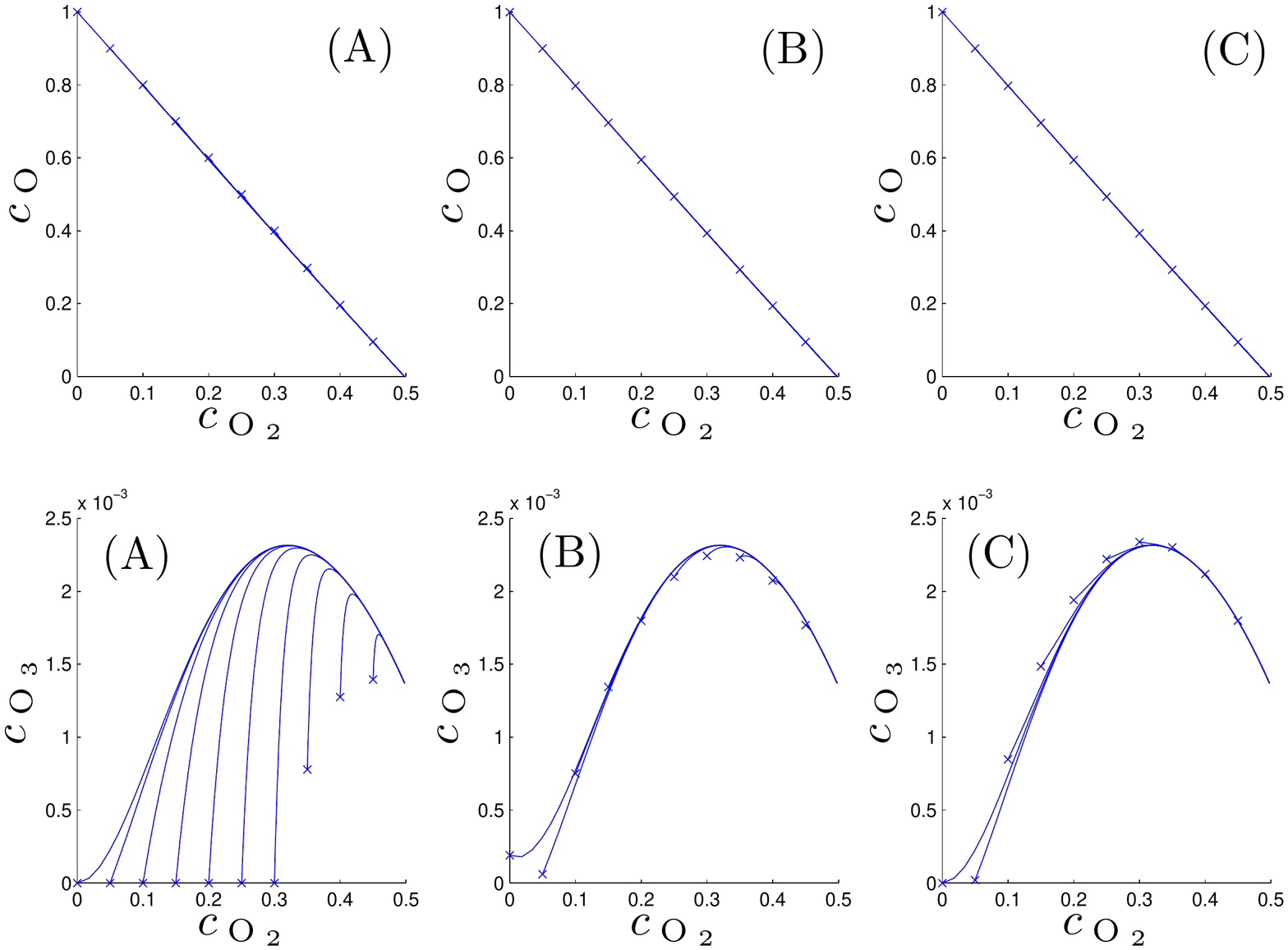}
  \end{center}
  \caption{\label{f:o3_35} Results for the ozone decomposition mechanism at a
    temperature of $T=350~{\mathrm{K}}$.}
\end{figure}

For $T=350~{\textrm{K}}$, the effects observed in Figure \ref{f:o3_50}
amplify. However, according to the results shown in Figure \ref{f:o3_35}
criteria B and C still perform reasonably well in this low-temperature region.

\subsubsection{Comparison with ILDM}
\begin{figure}
  \begin{center}
    \includegraphics[width=10.cm]{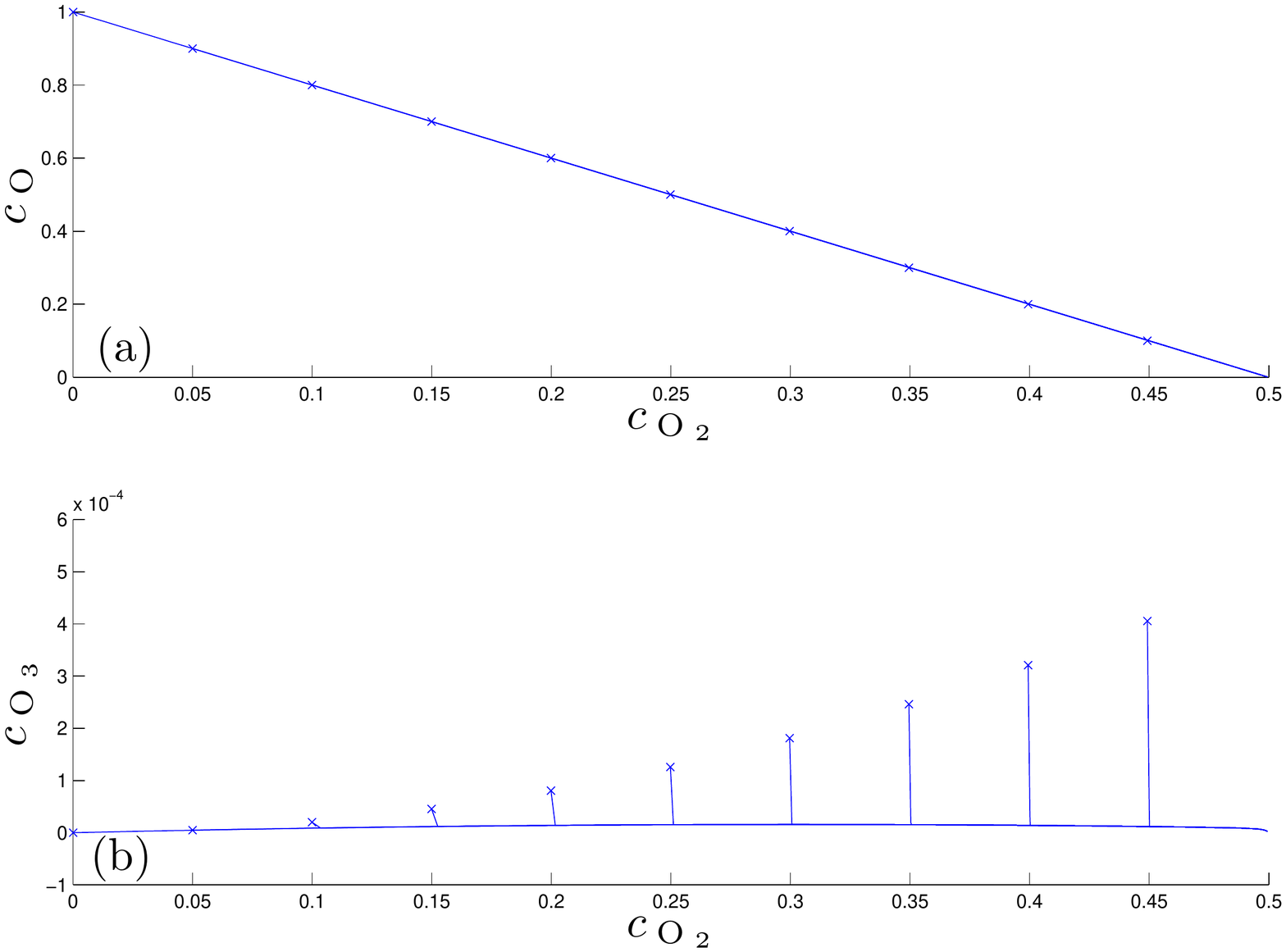}
  \end{center}
  \caption{\label{f:o3ildm} Results of the ILDM computation for the ozone
    decomposition mechanism at a temperature of $T=1000~{\mathrm{K}}$.  The
    x-marks depict the ILDM-points. They have been computed with the code used
    in \cite{Lebiedz2005c}, accuracy tolerance $10^{-9}$ for the solution of
    the ILDM-equation via Newton's method. The computation of ILDM points is
    initialized with the solution points of the optimization method using
    criterion B, cf.\ Fig.~\ref{f:o3_100}. Blue lines: trajectories started in
    the x-marks.}
\end{figure}

For the ozone decomposition mechanism we make a comparison of our results at
$T=1000~{\mathrm{K}}$ with the ILDM method \cite{Maas1992}. We numerically
compute ILDM-points for a range of $c_{\mathrm{O}_2}$ values. Figure
\ref{f:o3ildm} depicts the results. A comparison of Figure \ref{f:o3ildm}(b)
with the lower row of Figure \ref{f:o3_100} demonstrates a significantly
better performance of our trajectory optimization method. The ILDM points do
not lie close to the slow attracting manifold.

\subsubsection{Optimization Landscapes} \label{sss:o3_ol} We compute
optimization landscapes for the ozone decomposition model which has two
degrees of freedom. Initial values for $c_{\mathrm{O}_2}$ and
$c_{\mathrm{O}_3}$ are varied within the physically allowed range. The value
of the objective function is computed for trajectories corresponding to tuples
of initial values and depicted via color coding in a logarithmic scale. We
compare these optimization landscapes for $T=1000~{\mathrm{K}}$ and
$T=350~{\mathrm{K}}$.

\begin{figure}[htbp]
 \begin{center}
  \subfigure[ Optimization landscape for the ozone
    decomposition mechanism at $T=1000~{\mathrm{K}}$. The color represents the value of
    (\ref{eq:c2of}) for criterion A.]{\includegraphics[width=0.4\textwidth]{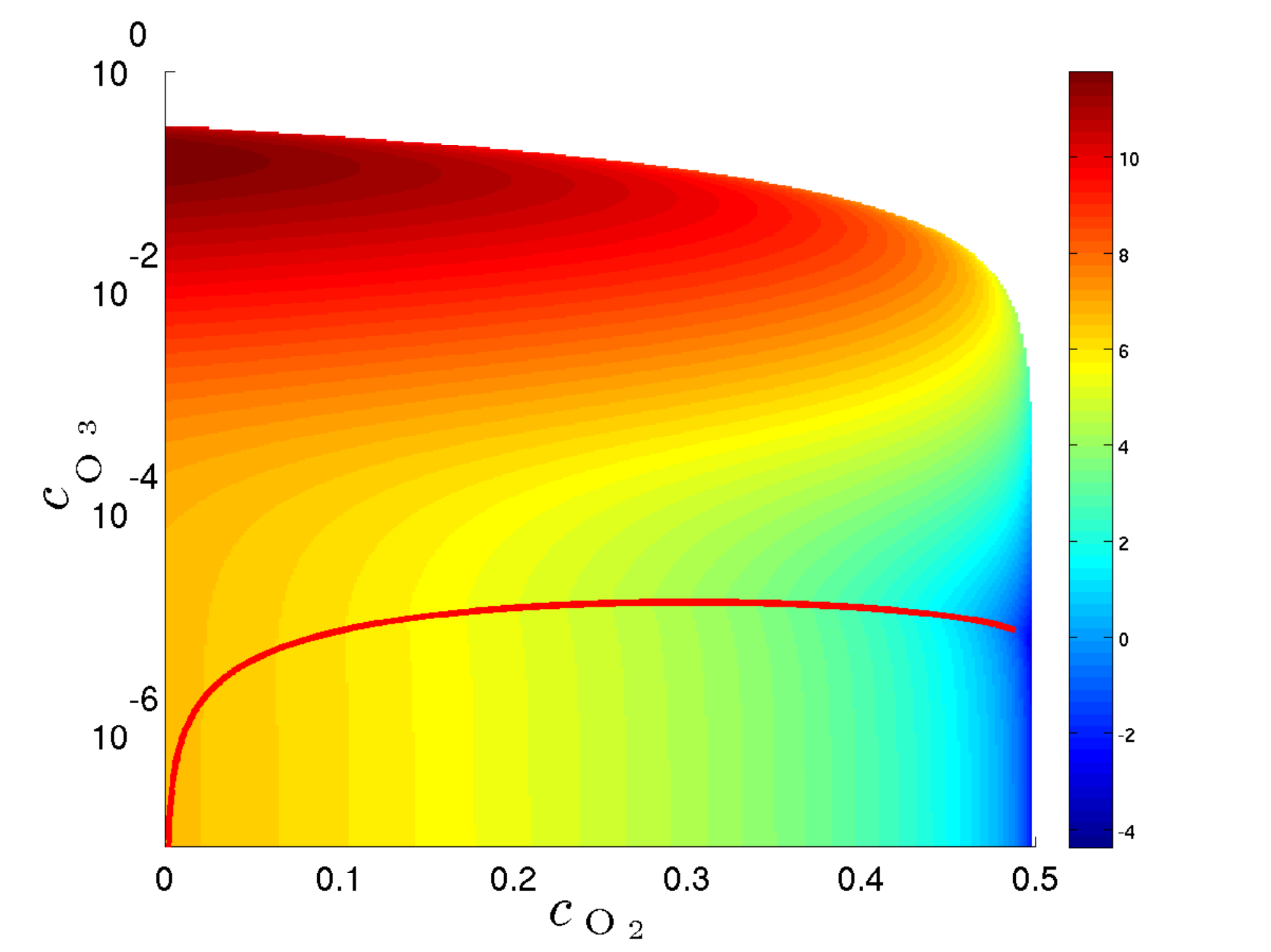}\label{f:o3_1000A}} \hspace*{0.1\textwidth}
  \subfigure[ Optimization landscape for the ozone mechanism
    at $T=1000~{\mathrm{K}}$: The color represents the value of
    (\ref{eq:cwof}) for criterion B.]{\includegraphics[width=0.4\textwidth]{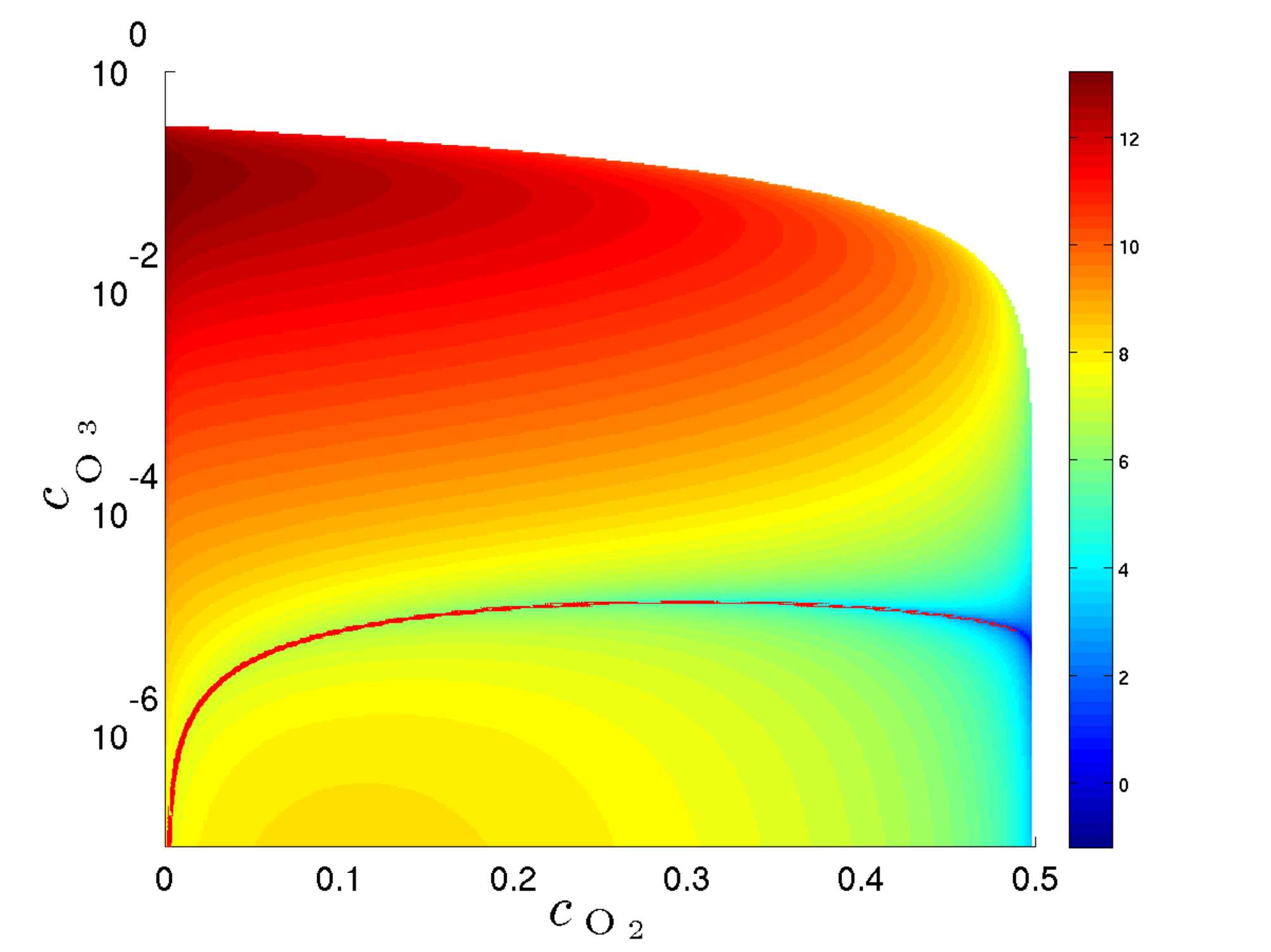}\label{f:o3_1000B}}
  \subfigure[ Optimization landscape for the ozone mechanism
    at $T=1000~{\mathrm{K}}$: The color represents the value of
    (\ref{eq:curv_tot}) for criterion C.]{\includegraphics[width=0.4\textwidth]{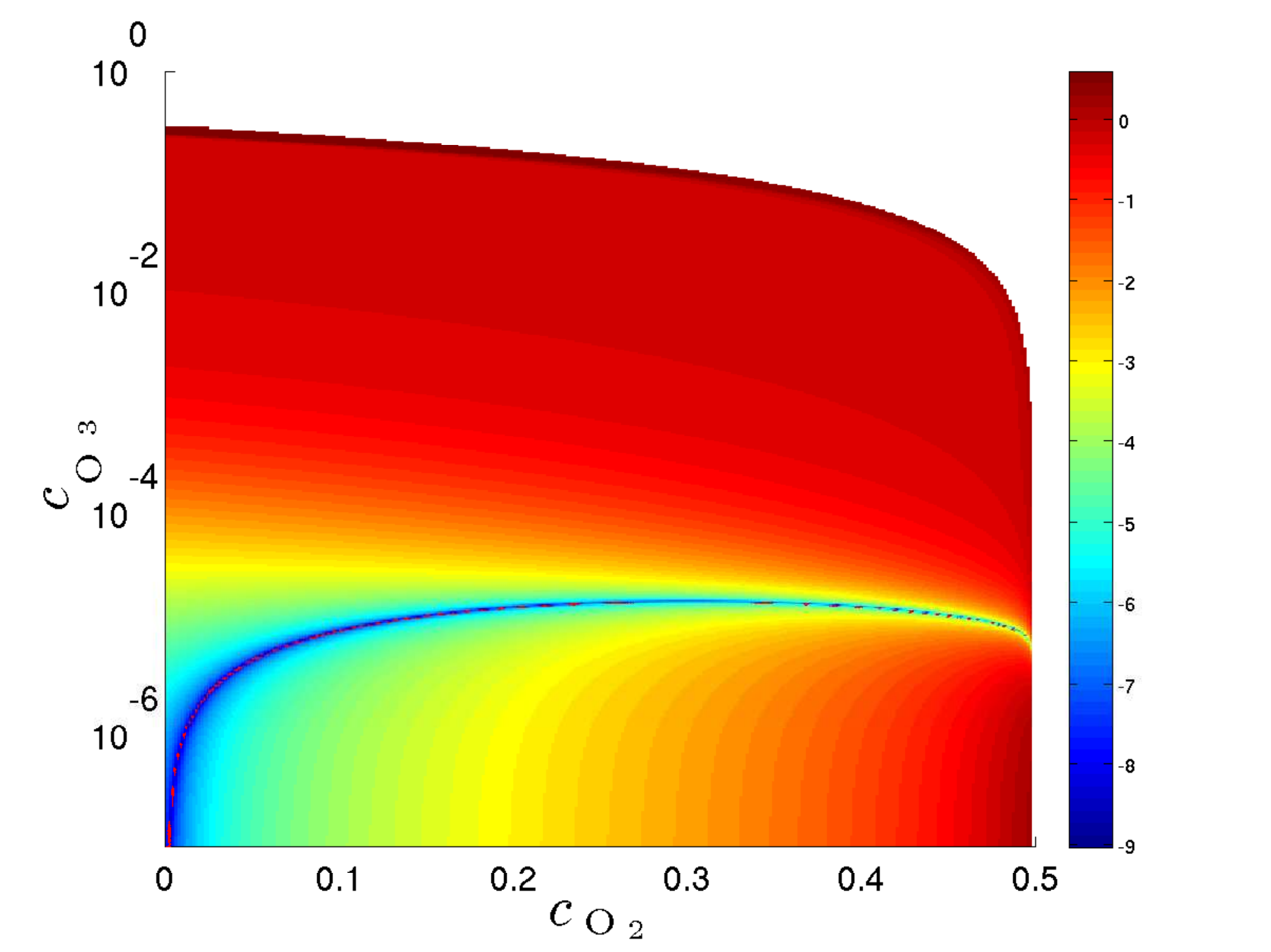}\label{f:o3_1000C}} \hspace*{0.1\textwidth}
  \subfigure[ Optimization landscape for the ozone mechanism at
    $T=350~{\mathrm{K}}$: criterion A.]{\includegraphics[width=0.4\textwidth]{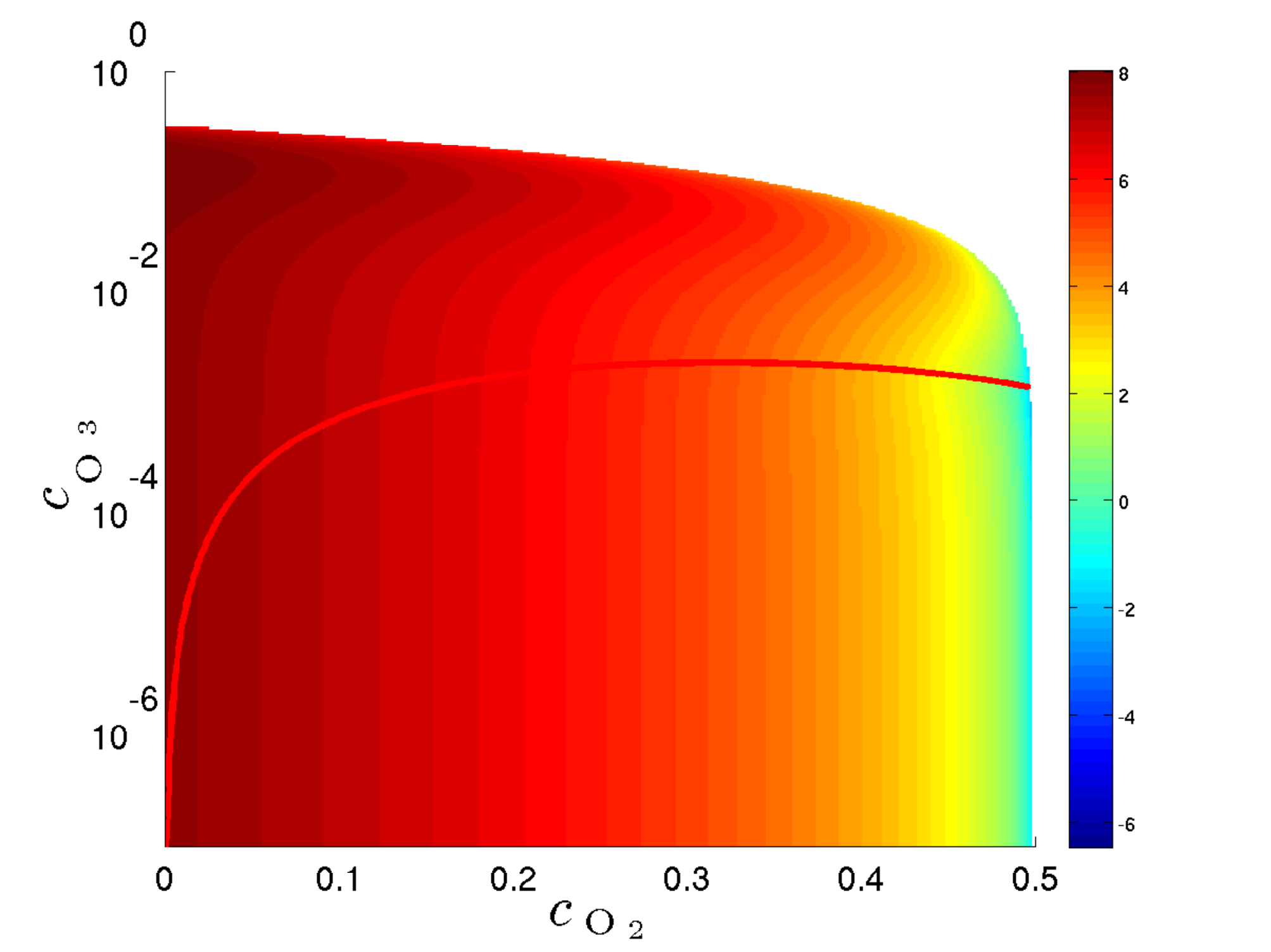} \label{f:o3_35A}}
  \subfigure[ Optimization landscape for the ozone mechanism at
    $T=350~{\mathrm{K}}$: criterion B.]{\includegraphics[width=0.4\textwidth]{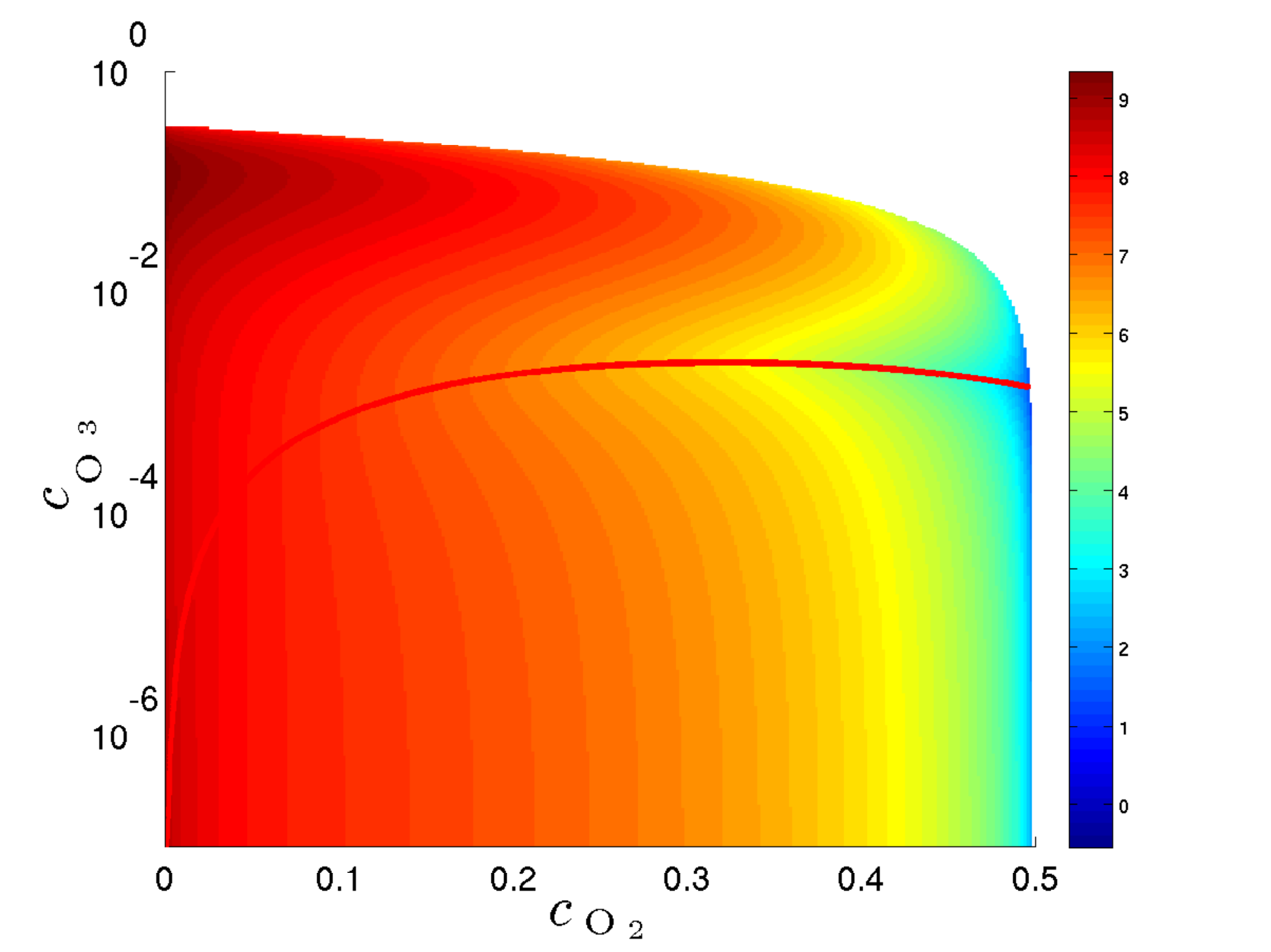} \label{f:o3_35B}} \hspace*{0.1\textwidth}
  \subfigure[ Optimization landscape for the ozone mechanism at
    $T=350~{\mathrm{K}}$: criterion C.]{\includegraphics[width=0.4\textwidth]{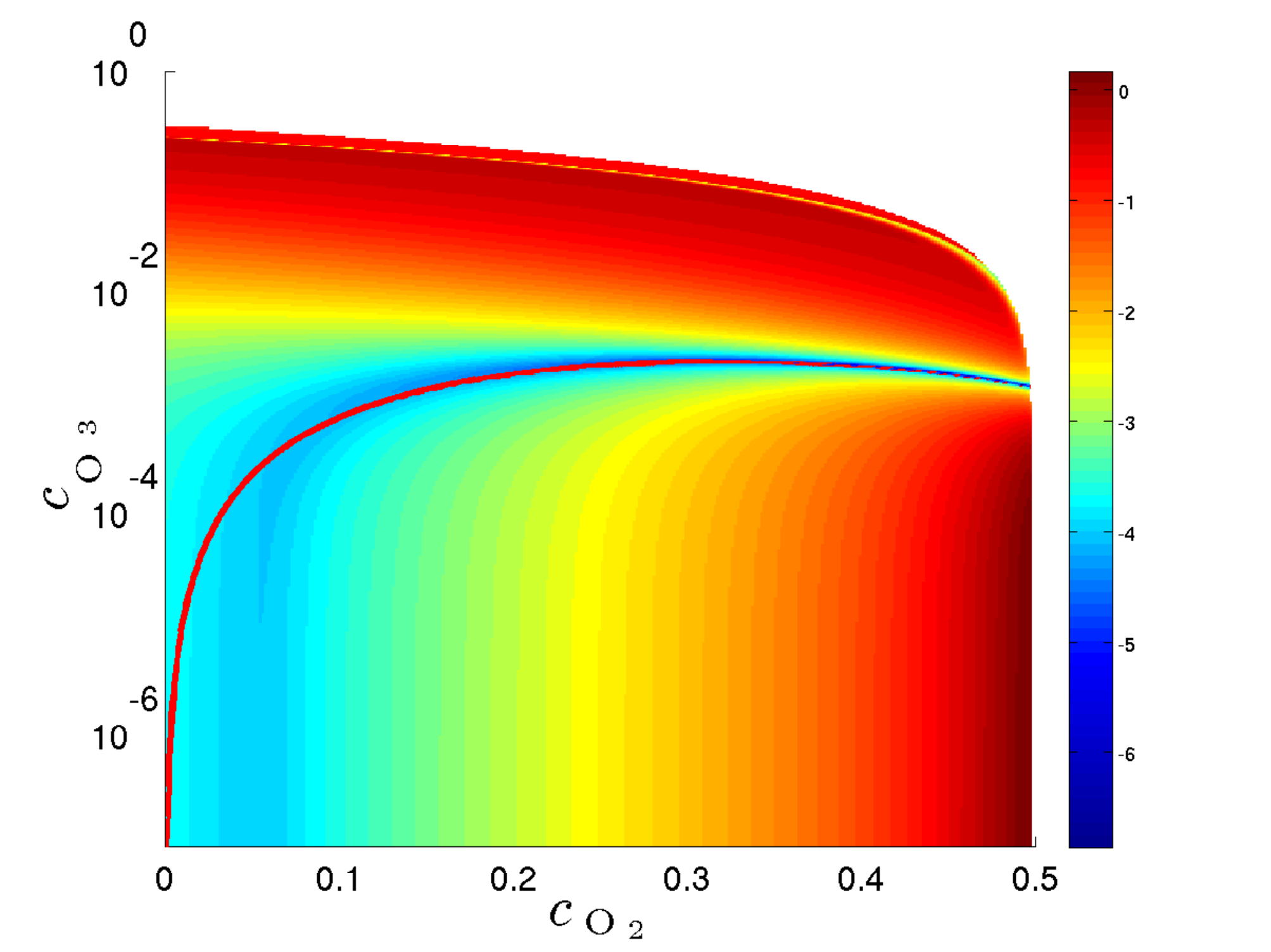} \label{f:o3_35C}}
 \end{center}
 \caption{ Optimization landscape for the ozone decomposition mechanism. The computed value of
    the different criteria for pairs of trajectory initial values is depicted in
    color with the numerical value corresponding to the logarithmically scaled
    colorbar. The value of\/ $c_{\mathrm{O}_3}$ is also logarithmically
    scaled. A trajectory (red) close to the SIM is shown.}
\end{figure}

Figure \ref{f:o3_1000A} shows the optimization landscape computed for
criterion A ($T=1000~{\mathrm{K}}$). The other criteria, B and C, give rise to
a much more distinct minimum of the objective function, cf.\ Figure
\ref{f:o3_1000B} and \ref{f:o3_1000C}. In the case of $T=350~{\mathrm{K}}$ criterion A fails, cf.\
Fig.~\ref{f:o3_35A}, no minimum near the SIM is found. The distinct minima for
criteria B and C become shallow but still allow for an optimal solution close
to the SIM. Figures \ref{f:o3_35B} and \ref{f:o3_35C} correspond to criteria B
and C respectively and visualize the results of the optimization problem
presented in Section~\ref{sss:o3or}.

\section{Summary and Discussion}
We present various geometrically motivated criteria for the numerical
computation of trajectories approximating slow (attracting) invariant
manifolds (SIM) in chemical reaction kinetics. The key idea of our approach is
to approximately span the SIM by trajectories being solutions of an
optimization problem for initial values of these trajectories. The objective
functional is supposed to characterize the extent of relaxation of chemical
forces (being minimal in the optimal solution) along a reaction trajectory on
its way towards equilibrium. Three different criteria are proposed and
motivated. Whereas the first two criteria use the directional derivative with
respect to its own direction of the tangent vector field of the kinetic ODE
system evaluated in a suitable norm, the third criterion uses a classical
differential geometric definition of curvature of trajectories regarded as
curves in $\mathbb{R}^n$.

These criteria are tested with three different chemical reaction mechanisms:
the Davis--Skodje problem \cite{Davis1999}, a six-species kinetic model for
hydrogen combustion, and a realistic ozone decomposition mechanism including
temperature dependence via Arrhenius kinetics. In all cases the quality of the
results are evaluated and compared.

Comparisons with the widely used ILDM-method \cite{Maas1992} show that our
method bears promise for improvements of slow manifold computations in
applications. Even though our optimization criteria do not guarantee invariant
manifolds in general, the solutions in our test examples are close to
invariance. It would be possible to compute invariant approximations of 1-D
manifolds by computing an optimal trajectory for reaction progress variable
values far from equilibrium und regard the resulting trajectory as a whole as
a SIM approximation. Then the manifold would be invariant by definition as a
trajectory of the ODE system. For the example of a kinetic model for the
temperature-dependent ozone decomposition it is demonstrated that our approach
also works reasonably well in low-temperature regions $T \leqslant 1000~K$
where the ILDM largely fails.

\section*{Acknowledgments}
This work was supported by the German Research Foundation (DFG) through the
Collaborative Research Center (SFB) 568. The authors also thank
Prof.~Dr.~Dr.~h.c.\ Hans~Georg~Bock (IWR, Heidelberg) for providing
the MUSCOD-II-package along with DAESOL and Prof.~Dr.~Moritz~Diehl
(K.U.\ Leuven/Belgium) for support with the initial value embedding
add-on for MUSCOD-II. We especially thank Andreas Potschka (IWR, Heidelberg)
for continuous support concerning MUSCOD-II and many helpful hints
related to derivative computation and problem formulation and Dr.~Mario~Mommer (IWR, Heidelberg) for various discussions.




\end{document}